%% file: hook.tex
\documentclass[a4paper]{article}
\usepackage{makeidx}
\usepackage{latexsym}
\usepackage{amsfonts}
\usepackage{amsmath}
\usepackage{amsthm}
\usepackage{amssymb}
\usepackage{amscd}
\usepackage[dvips]{graphicx}
\usepackage{cases}
\usepackage{color}
\usepackage{eepic}
\usepackage{setspace}
\usepackage{calrsfs}
\usepackage{authblk}
\usepackage[dvips]{hyperref}
\theoremstyle{definition}
\newtheorem{theorem}{Theorem}[section]
\newtheorem{prop}[theorem]{Proposition}
\newtheorem{lemma}[theorem]{Lemma}
\newtheorem{corollary}[theorem]{Corollary}
\newtheorem{definition}[theorem]{Definition}

\newtheorem{conjecture}[theorem]{Conjecture}

\newenvironment{demo}[1]{%
  \trivlist
  \item[\hskip\labelsep
        {\bf #1.}]
}{%
\hfill\qedsymbol
  \endtrivlist
}
 \makeatletter
    
    \@addtoreset{equation}{section}
  \makeatother
\renewcommand{\mathcal}{\mathrsfs}

\title{
$(q,t)$-hook formula for Birds and Banners
}
\author[1]{Masao ISHIKAWA}
\affil[1]{\small Department of Mathematics, Faculty of Education, University of the Ryukyus, Nishihara, Okinawa 901-0213, Japan,
{\tt ishikawa@edu.u-ryukyu.ac.jp}
}
\date{
}
\def\defterm#1{{\sl #1}\/}

\newcommand\rdots{\mathinner{\mkern1mu\raise0pt\vbox{\kern7pt\hbox{.}}
     \mkern2mu\raise4pt\hbox{.}\mkern2mu\raise8pt\hbox{.}\mkern1mu}}
\def\covered{\mathinner{\mkern1mu\raise0pt\vbox{\kern7pt\hbox{$<$}}
     \mkern-4mu\raise2pt\hbox{.}\mkern2mu}}
\def\covers{\mathinner{\mkern1mu\raise0pt\vbox{\kern7pt\hbox{$>$}}
     \mkern-12mu\raise2pt\hbox{.}\mkern8mu}}

\def\Z{{\mathbb{Z}}}
\def\N{{\mathbb{N}}}

\def\F{{\mathbb{F}}}

\def\tx{{\widetilde{x}}}
\def\ty{{\widetilde{y}}}
\def\tz{{\widetilde{z}}}
\def\hPhi{{\widehat{\Phi}}}
\def\tPhi{{\widetilde{\Phi}}}

\def\A{{\mathcal{A}}}
\def\CL{{\mathcal{C}}}
\begin{document}
%
%
%
\maketitle
\kern-15pt
\begin{abstract}
We study Okada's conjecture on $(q,t)$-hook formula of general $d$-complete posets.
Proctor classified $d$-complete posets into 15 irreducible ones.
We try to give a case-by-case proof of Okada's $(q,t)$-hook formula conjecture
using the symmetric functions.
Here we give a proof of the conjecture for birds and banners
in which we use Gasper's identity for VWP-series ${}_{12}W_{11}$.
\end{abstract}
%
%
%
\input hook01.tex
%
%
\input hook02.tex
%
%
\input hook03.tex

%
%
\input hook04.tex

%
%
%
%
%
%
%
%
%
%
%
%
%
%
%
%
%
%
%
\bibliographystyle{abbrv}
%
%
\bibliography{hook}
\label{sec:biblio}
\end{document}

%% file: hook01.tex
%
%
%
\section{Introduction and the main results}
The aim of this paper is to prove Okada's multivariate
hook formula conjecture for birds and banners,
i.e., Theorem~\ref{th:Okada-Birds-Banners}.
His conjecture is for general $d$-complete posets,
and here we give a partial proof for birds and banners only.
Proctor \cite{Pro1} has classified $d$-complete posets
into 15 irreducible classes. Okada 
\cite{Oka} has made his conjecture for general $d$-complete posets.
he has proven two cases in his paper
and we settle two cases in this paper 
so that the rest 11 classes are still left.
Even though we do a case-by-case proof,
we need the Macdonald polynomials and Gasper's identity for very well-poised series ${}_{12}W_{11}$ in our proof.
This paper is composed as follows.
In this section we recall the fundamental conceits
on $d$-complete partitions,
and then state our main result, i.e., 
Theorem~\ref{th:Okada-Birds-Banners}.
To sate Okada's conjecture we need the terminologies
on $d$-complete posets.
In Section~\ref{sec:Macdonald}
we recall the Macdonald polynomials.
In Section~\ref{sec:proof-theorems}
we rewrite Okada's conjecture by the Macdonald polynomials
and use the fact that the Macdonald polynomials
are basis for the ring of the symmetric functions.
In Section~\ref{sec:proof}
we prove the Macdonald polynomial identities obtained 
in Section~\ref{sec:proof-theorems}
using Gasper's identity.
\par\smallskip
%
%
%
%
%
%
Let $\N$ (resp. $\Z$) be the set of nonnegative integers (resp. integers).
Throughout this paper we use the standard  notation for $q$-series (see \cite{AAR,GR,KLS,KS}):
\begin{equation*}
    (a;q)_{\infty}=\prod_{k=0}^{\infty}(1-aq^{k}),\qquad
    (a;q)_{n}=\frac{(a;q)_{\infty}}{(aq^{n};q)_{\infty}}
\end{equation*}
for any integer $n$.
Usually  $(a;q)_{n}$ is called the  \defterm{$q$-shifted factorial},
and we frequently use the compact notation:
\begin{align*}
    &(a_{1},a_{2},\dots,a_{r};q)_{n}=(a_{1};q)_{n}(a_{2};q)_{n}\cdots(a_{r};q)_{n}.
\end{align*}
%
%
The \defterm{${}_{r+1}\phi_{r} $ basic hypergeometric series} is defined by
\begin{align}
{}_{r+1}\phi_{r}\left(\,
{{a_{1},a_{2},\dots,a_{r+1}}\atop{b_{1},\dots,b_{r}}};q,z
\,\right)
=\sum_{n=0}^{\infty}\frac{(a_{1},a_{2},\dots,a_{r+1};q)_{n}}{(q,b_{1},\dots,b_{r};q)_{n}}z^{n}.
\label{eq:phi-def}
\end{align}
%
%
A basic hypergeometric series ${}_{r+1}\phi_{r}$ is said to be 
\defterm{balanced} if
it satisfies $qa_{1}\cdots a_{r+1}=b_{1}\cdots b_{r}$ and $z=q$,
\defterm{well-poised} if
it satisfies $qa_{1}=a_{2}b_{1}=\cdots =a_{r+1} b_{r}$,
\defterm{very well-poised} if it is well-poised and satisfies
$b_{1}=a_{1}^{\frac12}$ and $b_{2}=-a_{1}^{\frac12}$
(see \cite[\S2.1]{GR}).
%
%
%
If ${}_{r+1}\phi_{r}$ is very well-poised series,
we use the notation
\begin{eqnarray*}
&{}_{r+1}W_{r}(a_{1};a_{4},\dots,a_{r+1};q,z)
={}_{r+1}\phi_{r}\Biggl[
{{a_{1},qa_{1}^{\frac12},-qa_{1}^{\frac12},a_{4},\dots,a_{r+1}}
\atop{a_{1}^{\frac12},-a_{1}^{\frac12},qa_{1}/a_{4},\dots,qa_{1}/a_{r+1}}}
\,;\,q,z
\Biggr].
\end{eqnarray*}
%
%
%
\begin{prop}
Gasper's identity (\cite[p.1065, (3.2)]{Gas}, \cite[pp.250, Ex.8.15]{GR}) reads as follows:
\begin{align}
&{}_{4}\phi_{3}\Biggl[
{{a,b,c,d}\atop{bq/a,cq/a,dq/a}}
;q,\frac{q^2}{a^2}
\Biggr]
=\frac
{(a/d.bq/d,cq/d,abc/d;q)_{\infty}}
{(q/d,ab/d,ac/d,bcq/d;q)_{\infty}}
\nonumber\\&\times
{}_{12}W_{11}\left(\frac{bc}{d};\left(\frac{bcq}{ad}\right)^{\frac12},-\left(\frac{bcq}{ad}\right)^{\frac12},
q\left(\frac{bc}{d}\right)^{\frac12},-q\left(\frac{bc}{d}\right)^{\frac12},
 \frac{ab}{d},\frac{ac}{d},a,b,c;q,\frac{q}{a}\right),
\label{eq:Gasper}
\end{align}
where at least one of $a$, $b$, $c$ is of the form $q^{-n}$ ($n=0,1,\dots$).
\end{prop}
We use the notation in \cite{Oka}.
For nonnegative integers $n$ and $m$ we write
\[
f(n;m)=f_{q,t}(n;m)=\frac{(t^{m+1};q)_{n}}{(t^{m};q)_{n}},
\]
and
\[
F(x)=F(x;q,t)=\frac{(tx;q)_{\infty}}{(x;q)_{\infty}},
\]
where $q$ and $t$ are parameters and $x$ is a variable (see \cite[(5)(6)]{Oka}).
Hereafter we use the convention that $f_{q,t}(n;m)=0$ for a negative integer $n<0$.
\par\bigbreak
%
%
%
%
We use the notation in \cite{Mac,Sta1} for partitions.
Let 
$\lambda= (\lambda_{1},\lambda_{2},\dots)$
 be a partition, 
i.e., 
$\lambda_{1}\geq\lambda_{2}\geq\dots$ with finitely many $\lambda_{i}$ unequal to zero.
The length and weight of $\lambda$, denoted by $\ell(\lambda)$ and $|\lambda|$, 
are the number and sum of the non-zero $\lambda_i$
respectively.
When $|\lambda|=N$ we say that $\lambda$ is
a partition of $N$, and the unique partition of zero is denoted by $\emptyset$. 
The multiplicity of the part $i$ in the partition $\lambda$ is denoted by $m_i(\lambda)$.
We identify a partition with its diagram (Ferrers graph)
\begin{equation}
D(\lambda)
=\{\,(i,j)\in\Z^{2}\,:\,1\leq j\leq\lambda_{i}\}.
\label{eq:Diagam}
\end{equation}
The conjugate $\lambda'$ of $\lambda$ is the partition obtained by
reflecting the diagram of $\lambda$ in the main diagonal.
A partition is said to be \defterm{strict} if we have strict inequalities 
$\lambda_{1}>\lambda_{2}>\cdots>\lambda_{r}> 0$ with $r=\ell(\lambda)$.
If $\lambda$ is a strict partition,
then its shifted diagram is defined by
\begin{equation}
S(\lambda)=
\{\,(i,j)\in\Z^{2}\,:\,i\leq j\leq \lambda_{i}+i-1\,\}.
\label{eq:ShiftedDiagam}
\end{equation}
Hereafter we may use the same symbol $\lambda$ to represent
its diagram (or shifted diagram).
\par
\begin{figure}[htbp]
\begin{center}
\setlength{\unitlength}{1.7pt}
\allinethickness{1pt}
\begin{picture}(30,90)
\put(15,5){\circle*{3}}
\put(15,15){\circle*{3}}
\put(15,35){\circle*{3}}
\put(5,45){\circle*{3}}
\put(25,45){\circle*{3}}
\put(15,55){\circle*{3}}
\put(15,75){\circle*{3}}
\put(15,85){\circle*{3}}
\put(15,5){\line(0,1){10}}
\put(15,15){\line(0,1){5}}
\multiput(15,20)(0,2){5}{\line(0,1){1}}
\put(15,30){\line(0,1){5}}
\put(15,35){\line(-1,1){10}}
\put(15,35){\line(1,1){10}}
\put(5,45){\line(1,1){10}}
\put(25,45){\line(-1,1){10}}
\put(15,55){\line(0,1){5}}
\multiput(15,60)(0,2){5}{\line(0,1){1}}
\put(15,70){\line(0,1){5}}
\put(15,75){\line(0,1){10}}
\put(10,70){\oval(3,30)[l]}
\put(-20,65){\makebox(30,10){$k-2$}}
\put(10,20){\oval(3,30)[l]}
\put(-20,15){\makebox(30,10){$k-2$}}
\put(15,80){\makebox(20,10){top}}
\put(-20,40){\makebox(25,10){side}}
\put(25,40){\makebox(25,10){side}}
\put(15,0){\makebox(35,10){bottom}}
\end{picture}
\end{center}
\caption{A double-tailed diamond poset $d_k(1)$
\label{fig:double-tailed}
}
\end{figure}
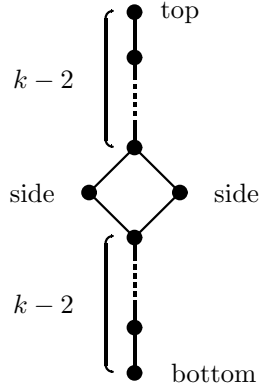

\bigbreak
We use standard notation and terminology of \cite[Chapter~3]{Sta1} related to posets.
We write $x\covered y$ if $x$ is covered by $y$,
i.e., $x<y$ and there is no $z\in P$ such that $x<z<y$.
A Hasse diagram is a diagram in which
one represents each element of $P$ as a vertex in the plane and draws an edge that goes upward from $x$ to $y$ whenever $y$ covers $x$.
\begin{definition}(\cite{Sta3}, \cite[\S3.15]{Sta1})
Let $P$ be a poset.
A \defterm{$P$-partition} is a map 
$\pi : P \to \N$ satisfying
\begin{equation}
x \le y \text{ in } P
\quad
\Longrightarrow
\quad
\pi(x) \ge \pi(y) \text{ in } \mathbb{N}.
\label{eq:P-partition}
\end{equation}
Let $\mathcal{A}(P)$ denote the set of $P$-partitions.
\end{definition}
%
%
%
%
First, we review the definition and some properties of d-complete posets. 
(See \cite{Pro1,Pro2}.) 
For $k\geq3$, 
we denote by $d_k(1)$ the poset consisting of $2k-2$ elements,
called \defterm{double-tailed diamond poset},
 with the Hasse diagram depicted in Figure~\ref{fig:double-tailed}.
The two incomparable elements are called the \defterm{sides},
the $k-2$ elements above them are called \defterm{neck} elements,
and the maximum and minimum elements are called \defterm{top} and \defterm{bottom} respectively.
If $k=3$ then we call $d_3(1)$ a \defterm{diamond}.
Let P be a poset. 
An interval $[w,v]=\{\,x\in P\,:\,w \leq x\leq v\,\}$
is called a \defterm{$d_k$-interval}
if it is isomorphic to $d_k(1)$
A \defterm{$d_k^-$-interval} ($k\geq4$) is an interval 
isomorphic to $d_k(1)-\{\text{top}\}$.
A \defterm{$d_3^-$-interval} consists of three elements 
$x$, $y$ and $w$ such that $w$ is covered by both $x$ and $y$.
A poset $P$ is \defterm{$d$-complete} if it satisfies
the following three conditions for every
$k\geq3$:
\begin{enumerate}
\item[(D1)]
If $I$ is a $d^{-}_{k}$-interval, 
then there exists an element $v$ such that $v$ covers the
maximal elements of $I$ and $I\cup\{v\}$ is a $d_{k}$-interval.
\item[(D2)]
If $I=[w,v]$ is a $d_k$-interval and the top $v$ covers $u$ in $P$, 
then $u\in I$.
\item[(D3)]
There are no $d^{-}_{k}$-intervals which differ only in the minimal elements.
\end{enumerate}
We quote a proposition due to Proctor \cite[Proposition in \S3]{Pro1} (also see \cite[Proposition~4.1]{Oka}):
%
%
%
%
\begin{prop}(\cite[Proposition in \S3]{Pro1})
\label{pr:Pro1} 
Let $P$ be a $d$-complete poset. 
Suppose that $P$ is connected,
i.e., the Hasse digram of $P$ is connected. 
Then we have
\begin{enumerate}
\item[(a)]
$P$ has a unique maximal element $v_{0}$.
\item[(b)]
For each $v\in P$, 
every saturated chain from $v$ to the maximum element $v_0$ has the same length. 
\end{enumerate}
Hence $P$ admits a rank function $r:P\to\N$ such that 
$r(x) =r(y) + 1$ if $x$ covers $y$.
\end{prop}
A \defterm{rooted tree} is a poset which has a unique maximal element,
and is such that each non-maximal element is covered by exactly one other element.
Let $P$ be a poset with a unique maximal element.
The top tree $T$ of $P$ is the filter 
(i.e., $x\in T$ and $y\geq x$ implies $y\in T$) of $P$, 
whose vertex set consists of all elements $x\in P$
such that every $y\geq x$ is covered by at most one other element of $P$.
$T$ is clearly a rooted tree and an element of $T$ is called \defterm{top tree element}.
Afterwards we use a particular kind of rooted tree.
Let $f\geq0$ and $h\geq g\geq0$ be integers.
The rooted tree $Y(f;g,h)$ consists of one “branch” element above which
a chain of $f$ elements has been adjoined and below which two non-adjacent chains with $g$ and $h$ elements, respectively.
\par\smallskip
Let $P$ be a connected $d$-complete poset with top tree $T$.
An element $x\in P$ is said to be \defterm{acyclic} 
if $x\in T$ and it is not in the neck of any $d_k$-interval for any $k\geq3$.
An element of $P$ is said to be \defterm{cyclic}
if it is not acyclic.
Let $Q$ be a $d$-complete poset containing an acyclic element $y$.
Let $P$ be a connected $d$-complete poset.
By Proposition~\ref{pr:Pro1} (a), 
let $x$ denote the unique maximal element of $P$.
Then the \defterm{slant sum} of $Q$ with $P$ at $y$, 
denoted $Q{}^{y}\backslash_{x}P$,
is the poset formed by creating a covering relation
$x\covered y$.
A $d$-complete poset $P$ is \defterm{slant irreducible} 
if it is connected and it cannot be expressed as a slant sum of two non-empty 
$d$-complete posets.
Suppose that $P$ is a connected $d$-complete poset with top tree $T$.
An edge $x\covered y$ of $P$ is a \defterm{slant edge} 
if $x,y\in T$ and $y$ is acyclic.
In \cite{Pro1}
Proctor proves $P$ is slant irreducible if and only if
 it contains no slant edges.
Also, $P$ is slant irreducible if and only if every acyclic
element is a minimal element of its top tree.
(\cite[Proposition~C of \S4]{Pro1})
Given  any connected $d$-complete poset $P$,
first locate all of its slant edges.
These may be erased in any order to produce a collection $P_1$, $P_2$,$\dots$
of uniquely determined smaller non-adjacent connected $d$-complete posets.
No new slant edges are created,
and so each of $P_1$, $P_2$,$\dots$ are slant irreducible.
We say that $P_1$, $P_2$,$\dots$ are the \defterm{slant irreducible components} of $P$.
If $P$ is an irreducible component,
then its top tree T is of the form $Y(f;g,h)$ for some $f\geq0$ and $h\geq g\geq1$ (\cite[Theorem of \S5]{Pro1}).
In the paper he establish the following theorem, 
which describe the structure of any connected $d$-complete poset.
%
%
%
%
\begin{prop}(Proctor \cite[Theorem in \S4]{Pro1})
Let $P$ be a connected $d$-complete poset.
It may be uniquely decomposed into a slant sum of one element posets 
and irreducible components.
The top tree of $P$ is an analogous slant sum of the top trees of the irreducible components.
\end{prop}
In \S7 of \cite{Pro1} Proctor defines 15 disjoint classes of irreducible components $\CL_{1},\dots,\CL_{15}$ and have shown that these 15 disjoint classes exhaust the set of all irreducible components.
For the list of 15 classes of irreducible $d$-complete posets
see \cite[Table~1]{Pro1}.
The diagram (\ref{eq:Diagam}) of an ordinary partition $\lambda$
or the shifted diagram (\ref{eq:ShiftedDiagam}) of a shifted partition $\lambda$
is regarded as a poset by defining its order structure as
\begin{equation}
(i_1,j_1)\geq(i_2,j_2)
\Longleftrightarrow
\text{ $i_1\leq i_2$ and $j_1\leq j_2$}.
\label{eq:order-shapes}
\end{equation}
By this order the poset represented by a diagram $P=D(\lambda)$ is called
a \defterm{shape} with its top tree $T=Y(f;g,h)$
where $f=0$, $g=\ell(\lambda)$ and $h=\ell(\lambda')$.
We use $\CL_1$ to express the class of shapes
which is a class of irreducible $d$-complete posets defined in \cite{Pro1}.
\par
\begin{figure}[htbp]
\begin{center}
\setlength{\unitlength}{1.0pt}
\begin{picture}(320,120)
\put( 10,110){\line(1,0){ 60}}
\put( 90,110){\line(1,0){ 80}}
\put(190,110){\line(1,0){ 80}}
\put(290,110){\line(1,0){ 20}}
\put( 30, 90){\line(1,0){ 40}}
\put( 90, 90){\line(1,0){ 80}}
\put(190, 90){\line(1,0){ 60}}
\put( 50, 70){\line(1,0){ 20}}
\put( 90, 70){\line(1,0){ 80}}
\put( 90, 50){\line(1,0){ 80}}
\put( 90, 30){\line(1,0){ 60}}
\put(110, 10){\line(1,0){ 20}}
\put( 75,110){\circle*{1}}
\put( 80,110){\circle*{1}}
\put( 85,110){\circle*{1}}
\put(175,110){\circle*{1}}
\put(180,110){\circle*{1}}
\put(185,110){\circle*{1}}
\put(275,110){\circle*{1}}
\put(280,110){\circle*{1}}
\put(285,110){\circle*{1}}
\put( 75, 90){\circle*{1}}
\put( 80, 90){\circle*{1}}
\put( 85, 90){\circle*{1}}
\put(175, 90){\circle*{1}}
\put(180, 90){\circle*{1}}
\put(185, 90){\circle*{1}}
\put( 75, 70){\circle*{1}}
\put( 80, 70){\circle*{1}}
\put( 85, 70){\circle*{1}}
\put(175, 70){\circle*{1}}
\put(180, 70){\circle*{1}}
\put(185, 70){\circle*{1}}
\put( 75, 50){\circle*{1}}
\put( 80, 50){\circle*{1}}
\put( 85, 50){\circle*{1}}
\put( 30, 90){\line(0,1){20}}
\put( 50, 70){\line(0,1){40}}
\put( 70, 50){\line(0,1){60}}
\put( 90, 50){\line(0,1){60}}
\put(110, 10){\line(0,1){20}}
\put(110, 50){\line(0,1){60}}
\put(130, 10){\line(0,1){20}}
\put(130, 50){\line(0,1){60}}
\put(150, 50){\line(0,1){60}}
\put(170, 50){\line(0,1){60}}
\put(190, 70){\line(0,1){40}}
\put(210, 90){\line(0,1){20}}
\put(230, 90){\line(0,1){20}}
\put(250, 90){\line(0,1){20}}
\put( 90, 35){\circle*{1}}
\put( 90, 40){\circle*{1}}
\put( 90, 45){\circle*{1}}
\put(110, 35){\circle*{1}}
\put(110, 40){\circle*{1}}
\put(110, 45){\circle*{1}}
\put(130, 35){\circle*{1}}
\put(130, 40){\circle*{1}}
\put(130, 45){\circle*{1}}
\put(150, 35){\circle*{1}}
\put(150, 40){\circle*{1}}
\put(150, 45){\circle*{1}}
\put( 10,110){\circle*{6}}
\put( 30,110){\circle*{6}}
\put( 50,110){\circle*{6}}
\put( 70,110){\circle*{6}}
\put( 90,110){\circle*{6}}
\put(110,110){\circle*{6}}
\put(130,110){\circle*{6}}
\put(150,110){\circle*{6}}
\put(170,110){\circle*{6}}
\put(190,110){\circle*{6}}
\put(210,110){\circle*{6}}
\put(230,110){\circle*{6}}
\put(250,110){\circle*{6}}
\put(270,110){\circle*{6}}
\put(290,110){\circle*{6}}
\put(310,110){\circle*{6}}
\put( 30, 90){\circle*{6}}
\put( 50, 90){\circle*{4}}
\put( 70, 90){\circle*{4}}
\put( 90, 90){\circle*{4}}
\put(110, 90){\circle*{4}}
\put(130, 90){\circle*{4}}
\put(150, 90){\circle*{4}}
\put(170, 90){\circle*{4}}
\put(190, 90){\circle*{4}}
\put(210, 90){\circle*{4}}
\put(230, 90){\circle*{4}}
\put(250, 90){\circle*{4}}
\put( 50, 70){\circle*{4}}
\put( 70, 70){\circle*{4}}
\put( 90, 70){\circle*{4}}
\put(110, 70){\circle*{4}}
\put(130, 70){\circle*{4}}
\put(150, 70){\circle*{4}}
\put(170, 70){\circle*{4}}
\put(190, 70){\circle*{4}}
\put( 70, 50){\circle*{4}}
\put( 90, 50){\circle*{4}}
\put(110, 50){\circle*{4}}
\put(130, 50){\circle*{4}}
\put(150, 50){\circle*{4}}
\put(170, 50){\circle*{4}}
\put( 90, 30){\circle*{4}}
\put(110, 30){\circle*{4}}
\put(130, 30){\circle*{4}}
\put(150, 30){\circle*{4}}
\put(110, 10){\circle*{4}}
\put(130, 10){\circle*{4}}
\put( 75, 45){\circle*{1}}
\put( 80, 40){\circle*{1}}
\put( 85, 35){\circle*{1}}
\end{picture}
\end{center}
\caption{Shifted shapes $C_{2}$
\label{fig:SShapes}
}
\end{figure}
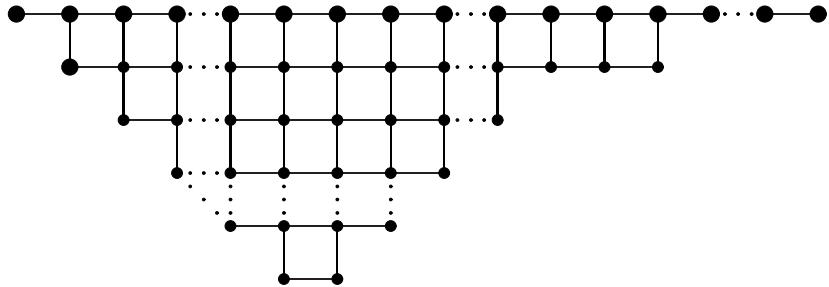

\smallskip
Another important class $\CL_2$ 
is the set of posets $P=S(\alpha)$ of shifted diagrams for strict partitions $\alpha$,
which is called \defterm{shifted shapes}
 with its top tree $T=Y(f,g,h)$
where $f=g=1$ and $h=\ell(\alpha)$.
Its Hasse diagram is designated by Figure~\ref{fig:SShapes}
in which the first row has $\alpha_1$ vertices,
the second row $\alpha_2$ vertices and so on.
When depicting these posets as a Hasse diagram,
we use the convention that a northwest vertex is
larger than another in southeast.
Here the larger dots and the heavier edges indicate the top tree.
For later use we denote by $P=P_{2}(\alpha)$ the Shifted shape associated
with a strict partition $\alpha$.
%
%
%
%
%
%
%
%
%
%
%
If $P=P_{2}(\alpha)$ is the shifted shape associated with a strict partition $\alpha$, then $P$-partition
\begin{equation}
\pi=(\pi_{ij})_{(i,j)\in S(\alpha)}
\label{eq:SShapes-partition}
\end{equation}
satisfies
\begin{equation}
\pi_{ij}\leq\pi_{i+1,j},
\qquad
\pi_{ij}\leq\pi_{i,j+1},
\label{eq:Shifted-PP-condition}
\end{equation}
whenever the both sides defined.
%
%
%
\begin{figure}[htbp]
\begin{center}
\setlength{\unitlength}{1.2pt}
\begin{picture}(160,80)
%
%
\put( 10, 70){\line(1,0){140}}
\put( 30, 50){\line(1,0){80}}
\put( 50, 30){\line(1,0){20}}
\put( 30, 50){\line(0,1){20}}
\put( 50, 30){\line(0,1){40}}
\put( 70, 10){\line(0,1){60}}
\put( 90, 50){\line(0,1){20}}
\put(110, 50){\line(0,1){20}}
\put( 10, 70){\circle*{6}}
\put( 30, 50){\circle*{6}}
\put( 30, 70){\circle*{6}}
\put( 50, 70){\circle*{6}}
\put( 70, 70){\circle*{6}}
\put( 90, 70){\circle*{6}}
\put(110, 70){\circle*{6}}
\put(130, 70){\circle*{6}}
\put(150, 70){\circle*{6}}
\put( 70, 10){\circle*{4}}
\put( 50, 30){\circle*{4}}
\put( 50, 50){\circle*{4}}
\put( 70, 30){\circle*{4}}
\put( 70, 50){\circle*{4}}
\put( 90, 50){\circle*{4}}
\put(110, 50){\circle*{4}}
\put( 13, 74){$\pi_{11}$}
\put( 33, 74){$\pi_{12}$}
\put( 53, 74){$\pi_{13}$}
\put( 73, 74){$\pi_{14}$}
\put( 93, 74){$\pi_{15}$}
\put(113, 74){$\pi_{16}$}
\put(133, 74){$\pi_{17}$}
\put(153, 74){$\pi_{18}$}
\put( 33, 54){$\pi_{22}$}
\put( 53, 54){$\pi_{23}$}
\put( 73, 54){$\pi_{24}$}
\put( 93, 54){$\pi_{25}$}
\put(113, 54){$\pi_{26}$}
\put( 53, 34){$\pi_{33}$}
\put( 73, 34){$\pi_{34}$}
\put( 73, 14){$\pi_{44}$}
\end{picture}
\end{center}
\caption{$P$-partition for shifted shape $(8,5,2,1)$
\label{fig:Shiftedshapes-partition}
}
\end{figure}
For example, Figure~\ref{fig:Shiftedshapes-partition} is a 
$P$-partition for shifted shape $(8,5,2,1)$.
\par\smallskip
In this paper we mainly consider only two classes, i.e.,
birds $\CL_3$ (Figure~\ref{fig:Birds}) and banners $\CL_6$ (Figure~\ref{fig:Banners}).
%
%
%
%
\begin{figure}[htbp]
\begin{center}
\setlength{\unitlength}{0.9pt}
\begin{picture}(260,160)
\put( 10,150){\line(1,0){40}}
\put( 70,150){\line(1,0){80}}
\put(170,150){\line(1,0){40}}
\put(230,150){\line(1,0){20}}
\put(110,130){\line(1,0){40}}
\put(170,130){\line(1,0){20}}
\put(110,110){\line(1,0){20}}
\put(110, 90){\line(1,0){20}}
\put(110, 70){\line(1,0){20}}
\put( 55,150){\circle*{1}}
\put( 60,150){\circle*{1}}
\put( 65,150){\circle*{1}}
\put(155,150){\circle*{1}}
\put(160,150){\circle*{1}}
\put(165,150){\circle*{1}}
\put(155,130){\circle*{1}}
\put(160,130){\circle*{1}}
\put(165,130){\circle*{1}}
\put(215,150){\circle*{1}}
\put(220,150){\circle*{1}}
\put(225,150){\circle*{1}}
\put(110, 10){\line(0,1){20}}
\put(110, 50){\line(0,1){40}}
\put(110,110){\line(0,1){40}}
\put(130, 70){\line(0,1){20}}
\put(130,110){\line(0,1){40}}
\put(150,130){\line(0,1){20}}
\put(170,130){\line(0,1){20}}
\put(190,130){\line(0,1){20}}
\put(110, 35){\circle*{1}}
\put(110, 40){\circle*{1}}
\put(110, 45){\circle*{1}}
\put(110, 95){\circle*{1}}
\put(110,100){\circle*{1}}
\put(110,105){\circle*{1}}
\put(130, 95){\circle*{1}}
\put(130,100){\circle*{1}}
\put(130,105){\circle*{1}}
\put(130,130){\line(1,-1){40}}
\put(190, 70){\line(1,-1){20}}
\put(175, 85){\circle*{1}}
\put(180, 80){\circle*{1}}
\put(185, 75){\circle*{1}}
\put( 10,150){\circle*{6}}
\put( 30,150){\circle*{6}}
\put( 50,150){\circle*{6}}
\put( 70,150){\circle*{6}}
\put( 90,150){\circle*{6}}
\put(110, 10){\circle*{6}}
\put(110, 30){\circle*{6}}
\put(110, 50){\circle*{6}}
\put(110, 70){\circle*{6}}
\put(110, 90){\circle*{6}}
\put(110,110){\circle*{6}}
\put(110,130){\circle*{6}}
\put(110,150){\circle*{6}}
\put(130, 70){\circle*{4}}
\put(130, 90){\circle*{4}}
\put(130,110){\circle*{4}}
\put(130,130){\circle*{4}}
\put(130,150){\circle*{6}}
\put(150,130){\circle*{4}}
\put(150,150){\circle*{6}}
\put(170,130){\circle*{4}}
\put(170,150){\circle*{6}}
\put(190,130){\circle*{4}}
\put(190,150){\circle*{6}}
\put(210,150){\circle*{6}}
\put(230,150){\circle*{6}}
\put(250,150){\circle*{6}}
\put(150,110){\circle*{4}}
\put(170, 90){\circle*{4}}
\put(190, 70){\circle*{4}}
\put(210, 50){\circle*{4}}
\put(114,154){$v$}
\put(134,134){$w$}
\end{picture}
\end{center}
\caption{Birds $C_{3}$
\label{fig:Birds}
}
\end{figure}
%
%
%
%
%
%
%
%
Let $\alpha=(\alpha_1,\alpha_2)$ and $\beta=(\beta_1,\beta_2)$
be strict partitions such that
$\alpha_1>\alpha_2>0$ and $\beta_1>\beta_2>0$.
Define the \defterm{bird} $P=P_{3}(\alpha,\beta;f)$ by
\[
P=P_\text{H}\cup P_\text{R}\cup P_\text{L}\cup P_\text{T}
\]
where
\begin{align*}
&
P_\text{H}=\{\,(1,j)\,:\,-f+1\leq j\leq 1\},
\\&
P_\text{R}=\{\,(i,j)\,:\,i\leq j\leq \alpha_{i}+i-1\ (i=1,2)\,\},
\\&
P_\text{L}=\{\,(i,j)\,:\,j\leq i\leq \beta_{j}+j-1\ (j=1,2)\,\},
\\&
P_\text{T}=\{\,(i,i)\,:\,2\leq i\leq f+2\}
\end{align*}
as a set and
we regard it as a poset
by defining its order structure (\ref{eq:order-shapes}) 
if and only if the both of $(i_{1},j_{1})$ and $(i_{2},j_{2})$ are in 
$P_\text{H}\cup P_\text{R}\cup P_\text{L}$ or in $P_\text{T}$
(see \cite[Table~1 and Figure~5.3]{Pro1}).
We call $P_\text{H}$ the \defterm{head}, $P_\text{T}$ the \defterm{tail},
$P_\text{R}$ (resp. $P_\text{L}$) the \defterm{right} (resp. \defterm{left}) \defterm{wing} of $P$.
The Hasse diagram of a bird is as in Figure~\ref{fig:Birds}.
Strictly speaking, we have to impose the condition
$\alpha_{1}=\alpha_{2}+1$ and $\beta_{1}=\beta_{2}+1$ to let $P$ be slant irreducible, but here we don't need this condition.
For example, the left-picture in Figure~\ref{fig:Birds-Banners} stands for $P=P_{3}((4,3),(4,2);2)$.
%
%
%
%
\begin{figure}[htbp]
\begin{center}
\setlength{\unitlength}{1.2pt}
\begin{picture}(340,80)
%
%
\put( 10, 70){\line(1,0){100}}
\put( 50, 50){\line(1,0){60}}
\put( 50, 30){\line(1,0){20}}
\put( 50, 10){\line(0,1){60}}
\put( 70, 30){\line(0,1){40}}
\put( 90, 50){\line(0,1){20}}
\put(110, 50){\line(0,1){20}}
\put( 70, 50){\line(1,-1){40}}
\put( 10, 70){\circle*{6}}
\put( 30, 70){\circle*{6}}
\put( 50, 70){\circle*{6}}
\put( 70, 70){\circle*{6}}
\put( 90, 70){\circle*{6}}
\put(110, 70){\circle*{6}}
\put( 50, 10){\circle*{6}}
\put( 50, 30){\circle*{6}}
\put( 50, 50){\circle*{6}}
\put( 50, 70){\circle*{6}}
\put( 70, 30){\circle*{4}}
\put( 70, 50){\circle*{4}}
\put( 90, 30){\circle*{4}}
\put( 90, 50){\circle*{4}}
\put(110, 10){\circle*{4}}
\put(110, 50){\circle*{4}}
\put( 53, 74){$v$}
\put( 73, 54){$w$}
\put( 33, 74){$v_1$}
\put( 13, 74){$v_2$}
\put( 93, 34){$w_1$}
\put(113, 14){$w_2$}
%
%
\put(150, 70){\line(1,0){180}}
\put(190, 50){\line(1,0){100}}
\put(210, 30){\line(1,0){ 40}}
\put(230, 10){\line(1,0){ 20}}
\put(190, 50){\line(0,1){20}}
\put(210, 10){\line(0,1){60}}
\put(230, 10){\line(0,1){60}}
\put(250, 10){\line(0,1){60}}
\put(270, 50){\line(0,1){20}}
\put(290, 50){\line(0,1){20}}
\put(150, 70){\circle*{6}}
\put(170, 70){\circle*{6}}
\put(190, 50){\circle*{6}}
\put(190, 70){\circle*{6}}
\put(210, 10){\circle*{4}}
\put(210, 30){\circle*{4}}
\put(210, 50){\circle*{4}}
\put(210, 70){\circle*{6}}
\put(230, 10){\circle*{4}}
\put(230, 30){\circle*{4}}
\put(230, 50){\circle*{4}}
\put(230, 70){\circle*{6}}
\put(250, 10){\circle*{4}}
\put(250, 30){\circle*{4}}
\put(250, 50){\circle*{4}}
\put(250, 70){\circle*{6}}
\put(270, 50){\circle*{4}}
\put(270, 70){\circle*{6}}
\put(290, 50){\circle*{4}}
\put(290, 70){\circle*{6}}
\put(310, 70){\circle*{6}}
\put(330, 70){\circle*{6}}
\put(194, 74){$v$}
\put(214, 54){$w$}
\put(174, 74){$v_1$}
\put(154, 74){$v_2$}
\put(214, 34){$w_1$}
\put(214, 14){$w_2$}
\end{picture}
\end{center}
\caption{Bird $P=P_{3}((4,3),(3,2);2)$ and banner $P=P_{6}((9,6,3,2);2)$
\label{fig:Birds-Banners}
}
\end{figure}
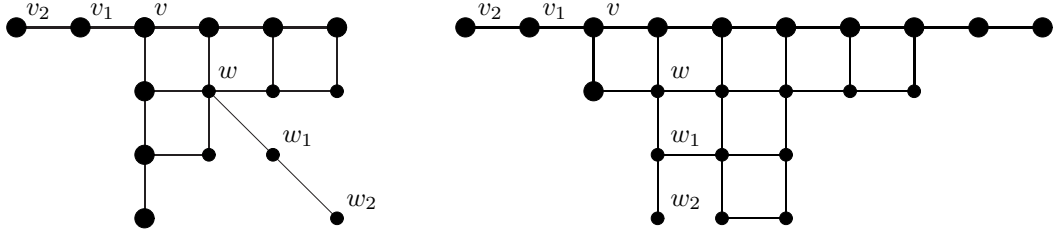
We have the chain $[v,v_2]$ (resp. $[w_2,w]$), which is the head (resp. tail) of $P$.
%
%
%
%
%
Recall that a $P$-partition $\pi$ satisfies the condition (\ref{eq:P-partition}). 
When $P=P_{3}(\alpha,\beta;f)$,
we associate the quadruple
$
(\sigma,\tau;\rho,\theta)
$
with $\pi$,
where
\[
\sigma=(\sigma_{i,j})_{(i,j)\in P_\text{R}},
\quad
\tau=(\tau_{i,j})_{(j,i)\in P_\text{L}},
\quad
\rho=(\rho_{i})_{i=0,\dots,f},
\quad
\theta=(\theta_{i})_{i=0,\dots,f}
\]
with
\begin{equation}
\begin{array}{llll}
\sigma_{i,j}=\pi(i,j)
&\text{ for $(i,j)\in P_\text{R}$,\qquad}
&\tau_{i,j}=\pi(j,i)
&\text{ for $(i,j)\in P_\text{L}$,}
\\
\rho_{-i+1}=\pi(1,i)
&\text{ for $(1,i)\in P_\text{H}$,\qquad}
&\theta_{i-2}=\pi(i,i)
&\text{ for $(i,i)\in P_\text{T}$.}
\end{array}
\label{eq:PP-condition-Birds}
\end{equation}
Hence we use the convention that $\rho_{0}=\sigma_{11}=\tau_{11}$
and $\theta_{0}=\sigma_{22}=\tau_{22}$.
We write $\pi=(\sigma,\tau;\rho,\theta)$ hereafter.
If $P=P_{3}((4,3),(4,2);2)$ then $\pi$ is as the left picture of Figure~\ref{fig:Birds-Banners-partition}.
\par\smallskip
%
%
%
%
%
%
%
Let $\alpha=(\alpha_1,\alpha_2,\alpha_3,\alpha_4)$ be a strict partition
such that $\alpha_1>\alpha_2>\alpha_3>\alpha_4>0$,
and let $f\geq2$ be a positive integer.
Let $P$ be the set
$P=P_\text{H}\cup P_\text{W}\cup P_\text{T}$
of lattice points in $\Z^2$,
where
\begin{align*}
&
P_\text{H}=\{\,(1,j)\,:\,-f+2\leq j\leq 1\},
\\&
P_\text{W}=\{\,(i,j)\,:\,i\leq j\leq \alpha_{i}+i-1\ (i=1,2,3,4)\,\},
\\&
P_\text{T}=\{\,(i,3)\,:\,3\leq i\leq f+2\}.
\end{align*}
We regard $P$ as a poset by defining the order relation (\ref{eq:order-shapes}) 
if both of  $(i_1,j_1)$ and $(i_2,j_2)$ are in $P_\text{H}\cup P_\text{W}$
or in $P_\text{T}$, and call it a \defterm{banner}
(see \cite[Table~1 and Figure~5.6]{Pro1}).
The Hasse diagram of a bird is as in Figure~\ref{fig:Birds} in general.
Strictly speaking again, we have to impose the condition
$\alpha_{1}=\alpha_{2}+1$ to let $P$ be slant irreducible, but we don't need this condition here.
We call $P_\text{H}$ the \defterm{head}, $P_\text{T}$ the \defterm{tail},
and $P_\text{W}$ the \defterm{wing} of $P$.
We use the symbol $P=P_{6}(\alpha;f)$ to mean the banner associated with 
a strict partition $\alpha$ and a positive integer $f$.
The Hasse diagram of a banner is given in Figure~\ref{fig:Banners}.
For example, the right picture in Figure~\ref{fig:Birds-Banners} stands for $P=P_{6}((9,6,3,2);2)$.
%
%
%
%
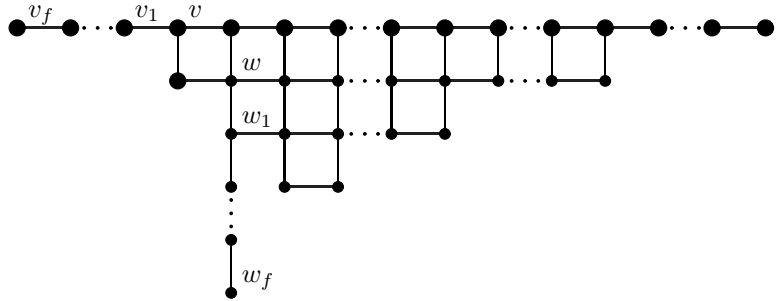
\begin{figure}[htbp]
\begin{center}
\setlength{\unitlength}{1.0pt}
\begin{picture}(300,120)
\put( 10,110){\line(1,0){20}}
\put( 50,110){\line(1,0){80}}
\put(150,110){\line(1,0){40}}
\put(210,110){\line(1,0){40}}
\put(270,110){\line(1,0){20}}
\put( 70, 90){\line(1,0){60}}
\put(150, 90){\line(1,0){40}}
\put(210, 90){\line(1,0){20}}
\put( 90, 70){\line(1,0){40}}
\put(110, 50){\line(1,0){20}}
\put(150, 70){\line(1,0){20}}
\put( 90,110){\line(1,0){20}}
\put( 35,110){\circle*{1}}
\put( 40,110){\circle*{1}}
\put( 45,110){\circle*{1}}
\put(135,110){\circle*{1}}
\put(140,110){\circle*{1}}
\put(145,110){\circle*{1}}
\put(195,110){\circle*{1}}
\put(200,110){\circle*{1}}
\put(205,110){\circle*{1}}
\put(255,110){\circle*{1}}
\put(260,110){\circle*{1}}
\put(265,110){\circle*{1}}
\put(135, 90){\circle*{1}}
\put(140, 90){\circle*{1}}
\put(145, 90){\circle*{1}}
\put(195, 90){\circle*{1}}
\put(200, 90){\circle*{1}}
\put(205, 90){\circle*{1}}
\put(135, 70){\circle*{1}}
\put(140, 70){\circle*{1}}
\put(145, 70){\circle*{1}}
\put( 70, 90){\line(0,1){20}}
\put( 90, 10){\line(0,1){20}}
\put( 90, 50){\line(0,1){60}}
\put(110, 50){\line(0,1){60}}
\put(130, 50){\line(0,1){60}}
\put(150, 70){\line(0,1){40}}
\put(170, 70){\line(0,1){40}}
\put(190, 90){\line(0,1){20}}
\put(210, 90){\line(0,1){20}}
\put(230, 90){\line(0,1){20}}
\put( 90, 35){\circle*{1}}
\put( 90, 40){\circle*{1}}
\put( 90, 45){\circle*{1}}
\put( 10,110){\circle*{6}}
\put( 30,110){\circle*{6}}
\put( 50,110){\circle*{6}}
\put( 70, 90){\circle*{6}}
\put( 70,110){\circle*{6}}
\put( 90, 10){\circle*{4}}
\put( 90, 30){\circle*{4}}
\put( 90, 50){\circle*{4}}
\put( 90, 70){\circle*{4}}
\put( 90, 90){\circle*{4}}
\put( 90,110){\circle*{6}}
\put(110, 50){\circle*{4}}
\put(110, 70){\circle*{4}}
\put(110, 90){\circle*{4}}
\put(110,110){\circle*{6}}
\put(130, 50){\circle*{4}}
\put(130, 70){\circle*{4}}
\put(130, 90){\circle*{4}}
\put(130,110){\circle*{6}}
\put(150, 90){\circle*{4}}
\put(150, 70){\circle*{4}}
\put(150,110){\circle*{6}}
\put(170, 70){\circle*{4}}
\put(170, 90){\circle*{4}}
\put(170,110){\circle*{6}}
\put(190, 90){\circle*{4}}
\put(190,110){\circle*{6}}
\put(210, 90){\circle*{4}}
\put(210,110){\circle*{6}}
\put(230, 90){\circle*{4}}
\put(230,110){\circle*{6}}
\put(250,110){\circle*{6}}
\put(270,110){\circle*{6}}
\put(290,110){\circle*{6}}
\put( 14,114){$v_f$}
\put( 54,114){$v_1$}
\put( 74,114){$v$}
\put( 94, 94){$w$}
\put( 94, 74){$w_1$}
\put( 94, 14){$w_f$}
\end{picture}
\end{center}
\caption{Banners $C_{6}$
\label{fig:Banners}
}
\end{figure}
%
%
%
%
%
%
%
If $P=P_{6}(\alpha;f)$ is the banner,
we associate a triplet $(\sigma;\rho,\theta)$ with a $P$-partition $\pi$, 
in which each component 
$\sigma=(\sigma_{i,j})_{(i,j)\in S(\alpha)}$,
$\rho=(\rho_{i})_{i=1,\dots,f}$, $\theta=(\theta_{i})_{i=1,\dots,f}$
are defined by
\begin{equation}
\begin{array}{llll}
\sigma_{i,j}=\pi(i,j)\quad
&\text{ for $(i,j)\in P_\text{W}$,}
\\
\rho_{i}=\pi(1,-i+2) 
&\text{ for $i=1,\dots,f$,}\quad
&\theta_{i}=\pi(i+2,3)\ 
&\text{ for $i=1,\dots,f$.}
\end{array}
\label{eq:PP-condition-Banners}
\end{equation}
Hence we have $\rho_{1}=\sigma_{11}$
and $\theta_{1}=\sigma_{33}$.
Hereafter we write $\pi=(\sigma;\rho,\theta)$.
For example, the right picture of Figure~\ref{fig:Birds-Banners-partition}
is a $P_{6}((9,6,3,2);2)$-partition.
%
%
%
%
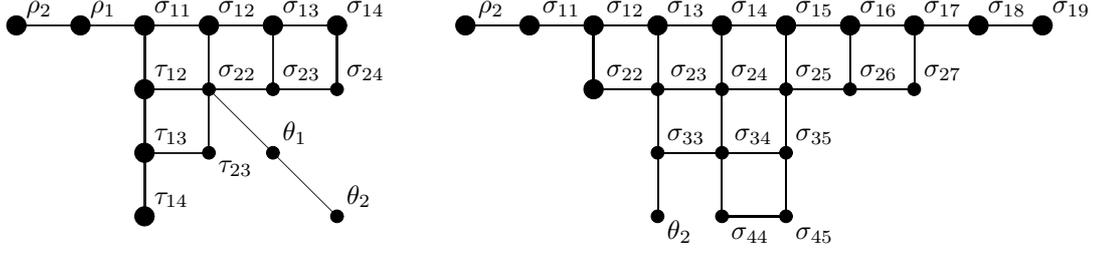
\begin{figure}[htbp]
\begin{center}
\setlength{\unitlength}{1.2pt}
\begin{picture}(340,80)
%
%
\put( 10, 70){\line(1,0){100}}
\put( 50, 50){\line(1,0){60}}
\put( 50, 30){\line(1,0){20}}
\put( 50, 10){\line(0,1){60}}
\put( 70, 30){\line(0,1){40}}
\put( 90, 50){\line(0,1){20}}
\put(110, 50){\line(0,1){20}}
\put( 70, 50){\line(1,-1){40}}
\put( 10, 70){\circle*{6}}
\put( 30, 70){\circle*{6}}
\put( 50, 70){\circle*{6}}
\put( 70, 70){\circle*{6}}
\put( 90, 70){\circle*{6}}
\put(110, 70){\circle*{6}}
\put( 50, 10){\circle*{6}}
\put( 50, 30){\circle*{6}}
\put( 50, 50){\circle*{6}}
\put( 50, 70){\circle*{6}}
\put( 70, 30){\circle*{4}}
\put( 70, 50){\circle*{4}}
\put( 90, 30){\circle*{4}}
\put( 90, 50){\circle*{4}}
\put(110, 10){\circle*{4}}
\put(110, 50){\circle*{4}}
\put( 53, 74){$\sigma_{11}$}
\put( 73, 54){$\sigma_{22}$}
\put( 33, 74){$\rho_1$}
\put( 13, 74){$\rho_2$}
\put( 93, 34){$\theta_{1}$}
\put(113, 14){$\theta_2$}
\put( 73, 74){$\sigma_{12}$}
\put( 93, 74){$\sigma_{13}$}
\put(113, 74){$\sigma_{14}$}
\put( 93, 54){$\sigma_{23}$}
\put(113, 54){$\sigma_{24}$}
\put( 53, 54){$\tau_{12}$}
\put( 53, 34){$\tau_{13}$}
\put( 53, 14){$\tau_{14}$}
\put( 73, 24){$\tau_{23}$}
%
%
\put(150, 70){\line(1,0){180}}
\put(190, 50){\line(1,0){100}}
\put(210, 30){\line(1,0){ 40}}
\put(230, 10){\line(1,0){ 20}}
\put(190, 50){\line(0,1){20}}
\put(210, 10){\line(0,1){60}}
\put(230, 10){\line(0,1){60}}
\put(250, 10){\line(0,1){60}}
\put(270, 50){\line(0,1){20}}
\put(290, 50){\line(0,1){20}}
\put(150, 70){\circle*{6}}
\put(170, 70){\circle*{6}}
\put(190, 50){\circle*{6}}
\put(190, 70){\circle*{6}}
\put(210, 10){\circle*{4}}
\put(210, 30){\circle*{4}}
\put(210, 50){\circle*{4}}
\put(210, 70){\circle*{6}}
\put(230, 10){\circle*{4}}
\put(230, 30){\circle*{4}}
\put(230, 50){\circle*{4}}
\put(230, 70){\circle*{6}}
\put(250, 10){\circle*{4}}
\put(250, 30){\circle*{4}}
\put(250, 50){\circle*{4}}
\put(250, 70){\circle*{6}}
\put(270, 50){\circle*{4}}
\put(270, 70){\circle*{6}}
\put(290, 50){\circle*{4}}
\put(290, 70){\circle*{6}}
\put(310, 70){\circle*{6}}
\put(330, 70){\circle*{6}}
\put(154, 74){$\rho_2$}
\put(174, 74){$\sigma_{11}$}
\put(194, 74){$\sigma_{12}$}
\put(213, 74){$\sigma_{13}$}
\put(233, 74){$\sigma_{14}$}
\put(253, 74){$\sigma_{15}$}
\put(273, 74){$\sigma_{16}$}
\put(293, 74){$\sigma_{17}$}
\put(313, 74){$\sigma_{18}$}
\put(333, 74){$\sigma_{19}$}
\put(194, 54){$\sigma_{22}$}
\put(214, 54){$\sigma_{23}$}
\put(233, 54){$\sigma_{24}$}
\put(253, 54){$\sigma_{25}$}
\put(273, 54){$\sigma_{26}$}
\put(293, 54){$\sigma_{27}$}
\put(213, 34){$\sigma_{33}$}
\put(234, 34){$\sigma_{34}$}
\put(253, 34){$\sigma_{35}$}
\put(213,  3){$\theta_{2}$}
\put(233,  3){$\sigma_{44}$}
\put(253,  3){$\sigma_{45}$}
\end{picture}
\end{center}
\caption{$P$-partitions
\label{fig:Birds-Banners-partition}
}
\end{figure}
\par\smallskip
Let $P$ be a connected $d$-complete poset and $T$ its top tree. 
Let $C$ be a set,
called \defterm{a set of colors},
whose cardinality is the same as $T$. 
A \defterm{coloring} of $P$ a coloring map $c$ of $P$
to the set of colors $C$.
$P$ is said to be \defterm{properly colored} if the coloring map $c$ satisfies
\begin{enumerate}
\item[(C1)]
$c(x)\neq c(y)$ if $x$ and $y$ are incomparable,
\item[(C2)]
$c(x)\neq c(y)$ if $x$ covers $y$.
\end{enumerate}
It is \defterm{simply colored} if, in addition: 
\begin{enumerate}
\item[(C3)]
whenever an interval $[w,v]$ is a chain,
the colors of the elements $c(x)$ in the interval $[w,v]$ are distinct.
\end{enumerate}
If $P$ is a rooted tree,
then it is simply colored
by the identity map $P\to P$,
i.e. we assign a distinct color to each vertex of $P$.
%
%
%
%
\begin{prop}(\cite[Proposition~8.6]{Pro2})
Let $P$ be a connected $d$-complete poset and $T$ its top tree.
Let $C$ be a set whose cardinality is the same as $T$. 
Then a bijection $c:T\to C$ can be uniquely extended 
to a proper coloring $c:P\to C$ satisfying the following condition:
\begin{enumerate}
\item[(C4)]
If $[w,v]$ is a $d_k$-interval then $c(w)=c(v)$.
\end{enumerate}
Such a map $c:P\to I$ is called a \defterm{$d$-complete coloring}.
\end{prop}
For example, in the both picture of Figure~\ref{fig:Birds-Banners}
because $[w_2,v_2]$ (resp. $[w_1,v_1]$, $[w,v]$) is a $d_5$-interval (resp.
$d_4$-interval, $d_3$-interval),
$w_2$ (resp. $w_1$, $w$) and $v_2$ (resp. $v_1$, $v$) have the same color.
In Figure~\ref{fig:Banners} $v_1$ (resp. $v_2$) and $v_3$ (resp.
$v_4$) have the same color since $[v_3,v_1]$ (resp. $[v_4,v_2]$ is a $d_4$-interval,
however,
the $v_1$ and $v_2$ have distinct colors since the both are in the top tree.
%
%
%
%
%
%
\begin{prop}
\label{pr:coloring} 
\begin{enumerate}
\item[(1)]
If $\alpha$ is a strict partition with length$\geq2$,
then the top tree of the shifted shape
$P=P_{2}(\alpha)$ is given by
\begin{equation}
T=\left\{(1,j)\,:\,1\leq j\leq\alpha_1\,\right\}\cup\{(2,2)\},
\end{equation}
and a $d$-complete coloring
$c:P\to\{0,0',1,2,\dots,\alpha_1-1\}$
is given by
\begin{equation}
c(i,j)=\begin{cases}
j-i&\text{ if $i<j$,}\\
0&\text{if $i=j$ and $i$ is odd,}\\
0'&\text{if $i=j$ and $i$ is even.}
\end{cases}
\end{equation}
Hence we see that $P$ has the top tree $Y(1;1,\alpha_1-1)$.
\item[(2)]
If $\alpha$ and $\beta$ are strict partitions with length$=2$
and $f\geq1$
then the top tree of the bird
$P=P_{3}(\alpha,\beta;f)$ is given by
\begin{equation}
T=\left\{(1,j)\,:\,-f+1\leq j\leq\alpha_1\,\right\}
\cup
\left\{(i,1)\,:\,1\leq i\leq\beta_1\,\right\},
\end{equation}
and a $d$-complete coloring
$c:P\to\{-f,\dots,-1,0,1,2,\dots,\alpha_1-1\}\cup\{1',2',\dots,(\beta_1-1)'\}$
is given by
\begin{equation}
c(i, j)=\begin{cases}
j-i&\text{if $i<j$, i.e., $(i,j)\in P_\text{R}$,}\\
(i-j)'&\text{if $1\leq j<i$, i.e., $(i,j)\in P_\text{L}$,}\\
j-1&\text{if $i=1$ and $j\leq1$, i.e., $(i,j)\in P_\text{H}$,}\\
-i+2&\text{if $i=j\geq2$, i.e., $(i,j)\in P_\text{T}$.}
\end{cases}
\end{equation}
Hence we see that 
$P$ has the top tree $Y(f;\alpha_1-1,\beta_1-1)$.
\item[(3)]
If $\alpha$ is a strict partitions with length$=4$
and $f\geq2$
then the top tree of the banner
$P=P_{3}(\alpha;f)$ is given by
\begin{equation}
T=\left\{(1,j)\,:\,-f+2\leq j\leq\alpha_1\,\right\}
\cup
\left\{(2,2)\right\},
\end{equation}
and a $d$-complete coloring
$c:P\to\{-f+1,\dots,-1,0,1,2,\dots,\alpha_1-1\}\cup\{0'\}$
is given by
\begin{equation}
c(i, j)=\begin{cases}
j-i&\text{ if $i\neq j$,}\\
0&\text{if $i=j$ and $i$ is odd,}\\
0'&\text{if $i=j$ and $i$ is even.}
\end{cases}
\end{equation}
Hence we see that 
$P$ has the top tree $Y(f;\alpha_1-1,1)$.
\end{enumerate}
\end{prop}
\begin{demo}{Proof}
(1) is obtained in \cite[Example~4.3(b)]{Oka}.
(2) and (3) are also obtained from (C1)--(C4).
\end{demo}
%
%
%
Let $P$ be a connected $d$-complete poset 
and $c:P\to C$ a $d$-complete coloring. 
Let $z_i$ ($i\in C$) be indeterminates.
For a $P$-partition $\pi\in\A(P)$, 
we write
\[
z^{\pi}=\prod_{v\in P}z^{\pi(v)}_{c(v)}.
\]
As in \cite[p.412]{Oka}
we associate a monomial
$z[H_{P}(v)]$ to each $v\in P$,
called the \defterm{hook monomial}, 
which is uniquely determined by induction as follows:
\begin{enumerate}
\item[(a)]
If $v$ is not the top of any $d_k$-interval, 
then we define
\[
z[H_{P}(v)]=\prod_{w\leq v}z_{c(w)}.
\]
\item[(b)] 
If $v$ is the top of a $d_k$-interval $[w,v]$, 
then we define
\[
z[H_{P}(v)]=\frac{z[H_{P}(x)]\cdot z[H_{P}(y)]}{z[H_{P}(w)]},
\]
where $x$ and $y$ are the sides of $[w,v]$.
\end{enumerate}
Further we denote $z[H_{p}]=\left\{\,z[H_{P}(v)]\,:\,v\in P\,\right\}$ the set of the hook monomials, and let $F\left(z[H_{p}];q,t\right)$ denote the product of $F\left(z[H_{P}(v)];q,t\right)$
 over $v\in P$, i.e.,
\[
F\left(z[H_{p}];q,t\right)=\prod_{v\in P}
F\left(z[H_{P}(v)];q,t\right). 
\]
\par\smallskip
Let $P$ be a connected $d$-complete poset with the maximum element $v_0$,
and the rank function $r:P\to\N$.
Let $T$ be the top tree of $P$. 
Take $T$ as a set of colors and let 
$c:P\to T$ be the $d$-complete coloring 
such that $c(v)=v$ for all $v\in T$.
Let $\widehat P=P\sqcup\{\widehat1\}$ be the extended poset,
where $\widehat1$ is the new maximum element of $\widehat P$ 
which covers $v_0$. 
Then $\widehat P$ has its top tree $\widehat T=T\sqcup\{\widehat1\}$,
where $\widehat c:\widehat P \to\widehat T$ with $\widehat c(\widehat1)=\widehat1$.
%
%
%
%
%
\begin{definition}
Given a $P$-partition $\pi\in\mathcal{A}(P)$,
let $\widehat\pi:\widehat P\to\N$ be the extensions of $\pi$ defined 
 by $\widehat\pi(\widehat1)=0$.
Define a weight $W_{P}(\sigma;q,t)$ by putting
\begin{align}
W_{P}(\pi;q,t)
=\frac{\displaystyle
\prod_{{x,y\in\widehat P}\atop{x<y,\ \widehat c(x)\sim \widehat c(y)}}
f(\pi(x)-\pi(y);d(x,y))
}
{\displaystyle
\prod_{{x,y\in P}\atop{x<y,\ c(x)=c(y)}}
f(\sigma(x)-\sigma(y);e(x, y))
f(\sigma(x)-\sigma(y);e(x, y)-1)
},
\label{eq:qt-weight}
\end{align}
where $\widehat c(x)\sim\widehat c(y)$ means that $\widehat c(x)$ and $\widehat c(y)$ are adjacent to each other in $T$, 
and
\begin{equation*}
d(x,y)=\frac{r(y)-r(x)-1}2,
\qquad
e(x,y)=\frac{r(y)-r(x)}2.
\end{equation*}
Note that if $c(x)\sim c(y)$ then $r(y)-r(x)$ is odd, 
and if $c(x) = c(y)$ then $r(y)-r(x)$ is even, 
hence $d(x,y)$ and $e(x,y)$ are nonnegative integers.
\end{definition}
Now we quote Okada's
$(q,t)$-hook formula conjecture. 
%
%
%
%
\begin{conjecture}(Okada \cite{Oka})
Let $P$ be a connected $d$-complete poset.
Using the notations defined
above, we have
\begin{align}
\sum_{\pi\in\A(P)}
W_{P}(\pi;q,t)z^{\pi}
=F\left(z[H_{p}];q,t\right).
\end{align}
\end{conjecture}
Okada has proven this conjecture for shapes and shifted shapes.
The purpose of this paper is to prove his conjecture for birds and banners.
%
%
%
%
\begin{theorem}
\label{th:Okada-Birds-Banners}
Okada's $(q,t)$-hook formula conjecture is true
for birds and banners.
\end{theorem}
%
%
%
%
Given a $P$-partition $\pi\in\A(P)$
for the shifted shape $P=P_{2}(\alpha)$ 
for a strict partition $\alpha$,
we write
\begin{align}
&f_{\alpha}^\text{ND}(\pi;q,t)
=\prod_{{(i,j)\in\alpha}\atop{i<j}}\prod_{m\geq0}
\frac
{f(\pi_{i,j}-\pi_{i-m,j-m-1};m)f(\pi_{i,j}-\pi_{i-m-1,j-m};m)}
{f(\pi_{i,j}-\pi_{i-m,j-m};m)f(\pi_{i,j}-\pi_{i-m-1,j-m-1};m)},
\label{eq:f-ND}
\\&
f_{\alpha}^\text{D}(\pi;q,t)
=\prod_{(i,i)\in\alpha}\prod_{{m\geq0}\atop{m\text{ even}}}
\frac
{f(\pi_{i,i}-\pi_{i-m-1,i-m};m)f(\pi_{i,i}-\pi_{i-m-2,i-m-1};m+1)}
{f(\pi_{i,i}-\pi_{i-m,i-m};m)f(\pi_{i,i}-\pi_{i-m-2,i-m-2};m+1)}.
\label{eq:f-D}
\end{align}
Here we use the convention that $\pi_{i,j}=0$ if $i\leq0$ or $j\leq0$.
Further we use the following short notation.
Let $m$ and $n$ be positive integers such that $m\leq n$.
When $\rho=(\rho_{m},\dots,\rho_{n})$ and $\theta=(\theta_{m},\dots,\theta_{n})$ 
satisfy
\begin{equation}
0\leq\rho_{n}\leq\cdots\leq\rho_{m}\leq\theta_{m}\leq\cdots\leq\theta_{n},
\label{eq:cond-rho-theata}
\end{equation}
we write
\begin{align}
\Phi_{m}^{n}(\rho,\theta;q,t)&=
\prod_{i=m+1}^{n}\frac
{f(\rho_{i-1}-\rho_{i};0)f(\theta_{i-1}-\rho_{i};0)f(\theta_{i}-\rho_{i-1};0)f(\theta_{i}-\theta_{i-1};0)}
{f(\theta_{i}-\rho_{i};i)f(\theta_{i}-\rho_{i};i+1)}.
\label{eq:def-Phi}
\end{align}
%
%
%
%
\begin{prop}
\label{pr:P-partitions}
\begin{enumerate}
\item[(1)]
Let $\alpha$ be a strict partition of length $r$ and $P=P_{2}(\alpha)$
the associated shifted shape.
If $\pi=(\pi_{ij})_{(i,j)\in\alpha}$ is a $P$-partition (\ref{eq:SShapes-partition}) satisfying the condition (\ref{eq:Shifted-PP-condition}),
then its weight $W_{P}(\pi;q,t)$ is given by
\begin{align}
W_{P}(\pi;q,t)
&=
f_{\alpha}^\text{D}(\pi;q,t)\,f_{\alpha}^\text{ND}(\pi;q,t).
\label{eq:weight-SShapes}
\end{align}
\item[(2)]
Let $\alpha$ 
and $\beta$
be strict partitions of length $2$.
Let $f>0$ be a positive integer,
and set $P=P_{3}(\alpha,\beta;f)$ to be the bird associated with
$\alpha$, $\beta$ and $f$.
If $\pi=(\sigma,\tau;\rho,\theta)$ is a $P$-partition 
 satisfying the condition (\ref{eq:PP-condition-Birds}),
then its weight $W_{P}(\pi;q,t)$ is given by
\begin{align}
W_{P}(\pi;q,t)
&=
\frac{f(\sigma_{22}-\sigma_{12};0)f(\tau_{22}-\tau_{12};0)
f(\rho_f;0)f(\theta_{f};f+1)}
{f(\sigma_{22}-\sigma_{11};0)f(\sigma_{22}-\sigma_{11};1)}
\nonumber\\&\times
\Phi_{0}^{f}(\rho,\theta;q,t)
f_{\alpha}^\text{ND}(\sigma;q,t)
f_{\beta}^\text{ND}(\tau;q,t).
\label{eq:weight-Birds}
\end{align}
Here we use the convention that 
$\sigma_{11}=\tau_{11}=\rho_{0}$ and $\sigma_{22}=\tau_{22}=\theta_{0}$.
\item[(3)]
Let $\alpha$ 
be a strict partition of length $4$.
Let $P=P_{6}(\alpha;f)$ be the banner associated with
$\alpha$ and $f$.
If $\pi=(\sigma;\rho,\theta)$ is a $P$-partition 
satisfying the condition (\ref{eq:PP-condition-Banners}),
then its weight $W_{P}(\pi;q,t)$ is given by
\begin{align}
W_{P}(\pi;q,t)
&=
f(\rho_f;0)f(\theta_{f};f+1)
\Phi_{1}^{f}(\rho,\theta;q,t)
f_{\alpha}^\text{D}(\sigma;q,t)
f_{\alpha}^\text{ND}(\sigma;q,t).
\label{eq:weight-Banners}
\end{align}
Here we use the convention that 
$\sigma_{11}=\rho_{1}$ and $\sigma_{33}=\theta_{1}$.
\end{enumerate}
\end{prop}
\begin{demo}{Proof}
For (1) (\ref{eq:weight-SShapes}) is exactly the same as \cite[Theorem~1.2 (9)]{Oka}.
For (2) and (3), one can compute $W_{P}(\pi;q,t)$ directly from the definition
(\ref{eq:qt-weight}).
\end{demo}
%
%
%
%
%
%
%
%
%
%
%
\begin{prop}
\label{pr:hook-length}
\begin{enumerate}
\item[(1)]
Let $\alpha$ be a strict partition of length $r$ and $P=P_{2}(\alpha)$
the associated shifted shape.
Let $n$ be an integer such that $n\geq\alpha_{1}$,
and let $\alpha^{c}$ be the strict partition formed by the complement of $\alpha$ in $[n]$, i.e.,
\[
\{\alpha_{1},\dots,\alpha_{r}\}\cup
\{\alpha_{1}^{c},\dots,\alpha_{n-r}^{c}\}=[n].
\]
We write $y_{0}=z_{0'}$ (see Proposition~\ref{pr:coloring} (1)) hereafter.
Then we have
\begin{align}
F\left(z[H_{p}];q,t\right)
&=\prod_{\alpha_{i}^{c}<\alpha_{j}}
F\left(\tz_{\alpha_{i}^{c}}^{-1}\tz_{\alpha_{j}};q.t\right)
\prod_{i}F\left(\tz_{\alpha_{i}};q,t\right)
\prod_{i<j}F\left(w\,\tz_{\alpha_{i}}\tz_{\alpha_{j}};q,t\right),
\label{eq:product-SShapes}
\end{align}
where 
$
\begin{cases}
\text{$w=y_{0}/z_{0}$ and $\tz_{i}=\prod_{k=0}^{i-1}z_{k}$ ($i=1,\dots,n$).}
&\text{ if $r$ is odd,}\\
\text{$w=z_{0}/y_{0}$ and $\tz_{i}=y_{0}\prod_{k=1}^{i-1}z_{k}$ ($i=1,\dots,n$).}
&\text{ if $r$ is even.}
\end{cases}
$
\item[(2)]
Let $\alpha=(\alpha_1,\alpha_2)$ 
and $\beta=(\beta_1,\beta_2)$
be strict partitions of length $2$.
Let $f>0$ be a positive integer,
and set $P=P_{3}(\alpha,\beta;f)$ the bird associated with
$f$, $\alpha$ and $\beta$.
Let $m,n$ be integers such that $m\geq\ell(\alpha)$ and $n\geq\ell(\beta)$,
and let $\alpha^{c}$ (resp. $\beta^{c}$) be the strict partition formed by the complement of $\alpha$ (resp. $\beta$) in $[m]$ (resp. $[n]$).
We write $y_{i}=z_{i'}$ for $i=1,\dots,\beta_1-1$ and
$x_{i}=z_{-i}$ for $i=1,\dots,f$.
Further we may write $x_{0}=y_{0}=z_{0}$.
 (See Proposition~\ref{pr:coloring} (2)).
Then we have
\begin{align}
F\left(z[H_{p}];q,t\right)
&=\prod_{\alpha_{i}^{c}<\alpha_{j}}
F\left(\tz_{\alpha_{i}^{c}}^{-1}\tz_{\alpha_{j}};q.t\right)
\prod_{\beta_{i}^{c}<\beta_{j}}
F\left(\ty_{\beta_{i}^{c}}^{-1}\ty_{\beta_{j}};q.t\right)
\prod_{i=1}^{f}F\left(\tx_{i};q.t\right)
\nonumber\\&\times
\prod_{i=1}^{f}F\left(\frac{\tx_{0}^{2}}{\tx_{i}}
\prod_{k,l=1}^{2}\ty_{l}\tz_{k};q.t\right)
\prod_{i,j=1}^{2}F\left(\tx_{0}\ty_{\beta_{j}}\tz_{\alpha_{i}};q,t\right)
\label{eq:product-Birds}
\end{align}
where $\tx_{i}=\prod_{k=i}^{f}x_{i}$ for $i=0,\dots,f$, $\ty_{i}=\prod_{k=1}^{i-1}y_{k}$ for $i=1,\dots,n$,
and $\tz_{i}=\prod_{k=1}^{i-1}z_{k}$ for $i=1,\dots,m$.
\item[(3)]
Let $\alpha=(\alpha_1,\alpha_2,\alpha_3,\alpha_4)$ 
be a strict partition of length $4$.
Let $P=P_{6}(f;\alpha)$ the banner associated with
$\alpha$ and $\beta$.
Let $n$ be an integer such that $n\geq 4=\ell(\alpha)$,
and let $\alpha^{c}$ be the strict partition formed by the complement of $\alpha$ in $[n]$.
We write $y_{0}=z_{0'}$ and
$x_{i}=z_{-i+1}$ for $i=2,\dots,f$
 (see Proposition~\ref{pr:coloring} (3)).
Hereafter we may use the convention that $x_{1}=z_{0}$.
Then we have
\begin{align}
F\left(z[H_{p}];q,t\right)
&=\prod_{\alpha_{i}^{c}<\alpha_{j}}
F\left(\tz_{\alpha_{i}^{c}}^{-1}\tz_{\alpha_{j}};q.t\right)
\prod_{i=2}^{f}F\left(\tx_{i};q.t\right)
\prod_{i=2}^{f}F\left(\frac{\tx_{2}^{\,2}}{\tx_{i}}\,w^2\prod_{i=1}^{4}\tz_{\alpha_{i}};q.t\right)
\nonumber\\&\times
\prod_{i=1}^{4}F\left(\tz_{\alpha_{i}};q,t\right)
\prod_{1\leq i<j\leq 4}F\left(\tx_{2}\,w\,\tz_{\alpha_{i}}\tz_{\alpha_{j}};q,t\right),
\label{eq:product-Banners}
\end{align}
where $w=\frac{z_{0}}{y_{0}}$, 
$\tx_{i}=\prod_{k=i}^{f}x_{i}$ for $i=1,\dots,f$
 and $\tz_{i}=y_{0}\prod_{k=1}^{i-1}z_{k}$ for $i=1,\dots,n$.
\end{enumerate}
\end{prop}
\begin{demo}{Proof}
\begin{enumerate}
\item[(1)]
If $P=P_{2}(\alpha)$, then we have
\[
z[H_{P}(i,j)]=\begin{cases}
\tz_{\alpha_{i}}\tz_{\alpha_{j+1}}
&\text{ if $i\leq j<r$,}\\
\tz_{\alpha_{i}}
&\text{ if $i\leq j=r$,}\\
\tz_{\alpha_{\alpha_{1}-j+1}^{c}}^{-1}\tz_{\alpha_{i}}
&\text{ if $i\leq r<j$,}
\end{cases}
\]
as stated in \cite[\S3, Proof Theorem~1.2]{Oka}.
\item[(2)]
If $P=P_{3}(\alpha,\beta;f)$, then we have
\[
z[H_{P}(i,j)]=\begin{cases}
\frac{\tx_{0}^{\,2}}{\tx_{-j+1}}\prod_{k,l=1}^{2}\ty_{k}\tz_{l}
&\text{ if $i=1$ and $-f+1\leq j\leq0$,}\\
\tx_{0}\ty_{\beta_{j}}\tz_{\alpha_{i}}
&\text{ if $1\leq i,j\leq2$,}\\
\tz_{\alpha_{\alpha_{1}-j+1}^{c}}^{-1}\tz_{\alpha_{i}}
&\text{ if $1\leq i\leq 2<j$,}\\
\ty_{\beta_{\beta_{1}-i+1}^{c}}^{-1}\ty_{\beta_{j}}
&\text{ if $1\leq j\leq 2<i$,}\\
\tx_{i-2}
&\text{ if $3\leq i=j\leq f+2$.}
\end{cases}
\]
\item[(3)]
If $P=P_{6}(\alpha;f)$, then we have
\[
z[H_{P}(i,j)]=\begin{cases}
\frac{\tx_{2}^{\,2}}{\tx_{-j+2}}\frac{z_{0}^2}{y_{0}^2}\prod_{k=1}^{4}\tz_{k}
&\text{ if $i=1$ and $-f+2\leq j\leq0$,}\\
\frac{z_{0}}{y_{0}}\tx_{2}\tz_{\alpha_{i}}\tz_{\alpha_{j+1}}
&\text{ if $1\leq i\leq j<4$,}\\
\tz_{\alpha_{i}}
&\text{ if $1\leq i\leq j=4$,}\\
\tz_{\alpha_{\alpha_{1}-j+1}^{c}}^{-1}\tz_{\alpha_{i}}
&\text{ if $1\leq i\leq 4<j$,}\\
\tx_{i-2}
&\text{ if $3< i\leq f+2$ and $j=3$.}
\end{cases}
\]
\end{enumerate}
\end{demo}
%
%
%
%
%
%
%
%
%
%
%
%
%
%
%
%
%
%
%
%
%
%
%
%
%
%
%
%
%
%
%
%
%
%
%
%
%
%
%
%
%
%
%
%
%

%% file: hook02.tex
%
%
%
%
%
%
%
%
%
%
\section{Macdonald polynomials}
\label{sec:Macdonald}
In this section
we recall the fundamental properties of Macdonald'polynomials
and consider it's application.
Especially Theorem~\ref{th:gMacMahon} and it's corollary will play
an important role in the next section.
\par\smallskip
We follow the notation and terminology of \cite{Mac} for the symmetric functions.
If $\lambda$ and $\mu$ are partitions then $\mu\subseteq\lambda$ 
if $\mu$ is contained in  $\lambda$,
i.e., $\mu_{i}\leq\lambda_{i}$ for all $i\geq1$. 
If $\mu\subseteq\lambda$ then the skew-diagram $\lambda/\mu$ denotes the set-theoretic
difference between $\lambda$ and $\mu$, i.e., 
those squares of $\lambda$ not contained in $\mu$. 
The skew diagram $\lambda/\mu$
is a vertical $r$-strip if $|\lambda-\mu|=|\lambda|-|\mu|=r$
 and if, for all $i\geq1$, 
$\lambda_{i}\geq\mu_{i}$ is at most one, 
i.e., each row of $\lambda-\mu$ contains at most one square. 
The set of all vertical $r$-strips is denoted by $\mathcal{V}_{r}$ and the set of all 
vertical strips by $\mathcal{V}=\biguplus_{r=0}^{\infty}\mathcal{V}_r$. 
The skew diagram $\lambda/\mu$ is a horizontal $r$-strip 
if $|\lambda-\mu| = r$ and if, for all $i\geq1$, $\lambda_{i}'-\mu_{i}'$ is at most one, 
i.e., each column of $\lambda-\mu$ contains at most one square.
For two partitions $\lambda$ and $\mu$,
we write $\lambda\succ\mu$
if $\lambda\supset\mu$ and $\lambda/\mu$ is a horizontal strip.
Note that  $\lambda/\mu$ is a horizontal strip
if and only if $\lambda_{1}\geq\mu_{1}\geq\lambda_{2}\geq\mu_{2}\geq\dots$.
The set of all horizontal $r$-strips is denoted by $\mathcal{H}_r$ 
and the set of all horizontal strips by $\mathcal{H}$.
Let $s = (i,j)$ be a square in the diagram of $\lambda$, 
and let $a(s)$ and $l(s)$ be the arm-length and leg-length of $s$, given by
\[
a(s)=\lambda_{i}-j,\qquad l(s)=\lambda_{j}'-i
\]
Then we define the rational functions
let
\[
b_{\lambda}(s)=b_{\lambda}(s;q,t)
:=\begin{cases}
\frac{1-q^{a(s)}t^{l(s)+1}}{1-q^{a(s)+1}t^{l(s)}},&\text{ if $s\in\lambda$,}\\
1,&\text{ otherwise,}
\end{cases}
\]
and \cite[(3.6)]{LSW} \cite[VI.7 (6.19), VI.7 Ex.4]{Mac}
\begin{align}
&b_{\lambda}(q,t)
:=\prod_{s\in\lambda}b_{\lambda}(s;q,t)
=\prod_{i\geq1}\prod_{m\geq0}\frac{f_{q,t}(\lambda_{i}-\lambda_{i+m+1};m)}{f_{q,t}(\lambda_{i}-\lambda_{i+m};m)},
\label{eq:b}
\\
&b^\text{el}_{\lambda}(q,t)
:=\prod_{{s\in\lambda}\atop{l(s)\text{ even}}}b_{\lambda}(s;q,t)
=\prod_{i\geq1}\prod_{{m\geq0}\atop{m\text{ even}}}\frac{f_{q,t}(\lambda_{i}-\lambda_{i+m+1};m)}{f_{q,t}(\lambda_{i}-\lambda_{i+m};m)},
\label{eq:b-el}
\\
&b^\text{oa}_{\lambda}(q,t)
:=\prod_{{s\in\lambda}\atop{a(s)\text{ odd}}}b_{\lambda}(s;q,t).
\end{align}
If $x=(x_1,x_2,\dots)$ and $y=(y_1,y_2,\dots)$ are two sequences of independent indeterminates,
then we write
\begin{equation}
\Pi(x;y;q,t)=\prod_{i,j}\frac{(tx_iy_j;q)_{\infty}}{(x_iy_j;q)_{\infty}}
=\prod_{i,j}F(x_iy_j;q,t).
\end{equation}
\par\bigbreak
Let $\mathfrak{S}_{n}$ denote the symmetric group, 
acting on $x = (x_{1},\dots, x_{n})$ by permuting the $x_i$,
and let $\Lambda_{n}=\mathbb{Z}[x_{1},\dots,x_{n}]^{\mathfrak{S}_{n}}$ and 
$\Lambda$ denote the ring of symmetric polynomials in $n$ independent
variables and the ring of symmetric polynomials in countably many variables, 
respectively.
For $\lambda=(\lambda_{1},\dots,\lambda_{n})$ a partition of at most $n$ parts 
the monomial symmetric function $m_{\lambda}$ is defined as
\[
m_{\lambda}(x)=\sum_{\alpha}x^{\alpha}
\]
where the sum is over all distinct permutations $\alpha$ of $\lambda$, 
and $x=(x_{1},\dots,x_{n})$.
For $\ell(\lambda)>n$
we set $m_{\lambda}(x)=0$. 
The monomial symmetric functions $m_{\lambda}(x)$ for $\ell(\lambda)\leq n$
 form a $\mathbb{Z}$-basis of $\Lambda_{n}$.
For $r$ a nonnegative integer the power sums $p_{r}$ are given by $p_0=1$ and 
$p_r=m_{(r)}$ for $r>1$.
More generally the power-sum products are defined as
$p_{\lambda}(x)=p_{\lambda_1}(x)p_{\lambda_2}(x)\cdots$ 
for an arbitrary partition $\lambda= (\lambda_{1},\lambda_{2},\dots)$.
Define the Macdonald scalar product $\langle\cdot,\cdot\rangle_{q,t}$ 
on the ring of symmetric functions by
\[
\langle p_{\lambda},p_{\mu}\rangle_{q,t}=\delta_{\lambda\mu}z_{\lambda}\prod_{i}
\prod^{n}_{i=1}
\frac{1-q^{\lambda_{i}}}{1-t^{\lambda_{i}}}
\]
with 
$
z_{\lambda}=\prod_{i\geq1}i^{m_{i}} m_{i}!
$
 and $m_{i}=m_{i}(\lambda)$.
If we denote the ring of symmetric functions in $\Lambda_{n}$
variables over the field $\mathbb{F}=\mathbb{Q}(q,t)$ of rational functions 
in $q$ and $t$ by $\Lambda_{n,\mathbb{F}}$, then the Macdonald
polynomial $P_{\lambda}(x)=P_{\lambda}(x;q,t)$ is the unique symmetric polynomial in $\Lambda_{n,\mathbb{F}}$
 such that [VI (4.7)]{Mac}:
%
%
%
\[
P_{\lambda}=\sum_{\mu\leq\lambda}u_{\lambda\mu}(q,t)m_{\mu}(x)
\]
with $u_{\lambda\lambda}=1$ and
\[
\langle P_{\lambda},P_{\mu}\rangle_{q,t}=0\qquad\text{ if $\lambda\neq\mu$.}
\]
The Macdonald polynomials $P_{\lambda}(x;q,t)$ with $\ell(\lambda)\leq n$ 
form an $\mathbb{F}$-basis of $\Lambda_{n,\mathbb{F}}$.
If $\ell(\lambda)>n$ then $P_{\lambda}(x;q,t)=0$.
$P_{\lambda}(x;q,t)$ is called \defterm{Macdonald's $P$-function}.
Since 
$P_{\lambda}(x_{1},\dots, x_{n},0;q,t)= P_{\lambda}(x_{1},\dots, x_{n};q,t)$ 
one can extend the Macdonald polynomials
to symmetric functions containing an infinite number of independent variables 
$x=(x_{1},x_{2},\dots)$,
to obtain a basis of $\mathbb{F}=\Lambda\otimes\mathbb{F}$.
A second Macdonald symmetric function, called \defterm{Macdonald's $Q$-function}, is defined as
\begin{equation}
Q_{\lambda}(x;q,t)=b_{\lambda}(q,t)P_{\lambda}(x;q,t).
\label{eq:Q=bP}
\end{equation}
The normalization of the Macdonald inner product is then
$\langle P_{\lambda},Q_{\mu}\rangle_{q,t}=\delta_{\lambda\mu}$ for
all $\lambda,\mu$,
which is equivalent to
\begin{equation}
\sum_{\lambda}
P_{\lambda}(x;q,t)Q_{\lambda}(y;q,t)=\Pi(x;y;q,t).
\label{eq:normalisation}
\end{equation}
(See \cite[VI.4, (4.13)]{Mac}.)
%
%
%
Let $g_{r}(x;q,t):= Q_{(r)}(x;q,t)$, or equivalently, \cite[VI.2, (2.8)]{Mac}
\[
\prod_{i=1}^{\infty}\frac{(tx_{i}y;q)_{\infty}}{(x_{i}y;q)_{\infty}}
=\sum_{r=0}^{\infty}g_{r}(x;q,t)y^{r}.
\]
Then the \defterm{Pieri coefficients} $\phi_{\lambda/\mu}$ and $\psi_{\lambda/\mu}$ 
are given by \cite[VI.6, (6.24)]{Mac}
\begin{align*}
&
P_{\mu}(x;q,t)g_{r}(x;q,t)
=
\sum_{{\lambda}\atop{\lambda-\mu\in\mathcal{H}_r}}
\phi_{\lambda/\mu}(q,t)P_{\lambda}(x;q,t),
\\&
Q_{\mu}(x;q,t)g_{r}(x;q,t)
=
\sum_{{\lambda}\atop{\lambda-\mu\in\mathcal{H}_r}}
\psi_{\lambda/\mu}(q,t)Q_{\lambda}(x;q,t).
\end{align*}
Another direct expressions for $\phi_{\lambda/\mu}$ and $\psi_{\lambda/\mu}$ is given
in \cite[VI.6, Ex.2]{Mac} as 
\begin{align}
&
\phi_{\lambda/\mu}(q,t)
=\prod_{1\leq i\leq j\leq\ell(\lambda)}
\frac{f(\lambda_{i}-\mu_{j};j-i)f(\mu_{i}-\lambda_{j+1};j-i)}
{f(\lambda_{i}-\lambda_{j};j-i)f(\mu_{i}-\mu_{j+1};j-i)},
\label{eq:Pierri-phi}
\\&
\psi_{\lambda/\mu}(q,t)
=\prod_{1\leq i\leq j\leq\ell(\mu)}
\frac{f(\lambda_{i}-\mu_{j};j-i)f(\mu_{i}-\lambda_{j+1};j-i)}
{f(\mu_{i}-\mu_{j};j-i)f(\lambda_{i}-\lambda_{j+1};j-i)}.
\label{eq:Pierri-psi}
\end{align}
Here we use these expressions to rewrite Okada's $(q,t)$-hook formula conjectures
by the Pieri coefficients.
For any three partitions $\lambda$, $\mu$, $\nu$
let $f^{\lambda}_{\mu\nu}$ be the coefficient $P_{\lambda}$
in the product $P_{\mu}P_{\nu}$: \cite[VI (7.1')]{Mac}:
\begin{equation}
P_{\mu}(x;q,t)P_{\mu}(x;q,t)=\sum_{\lambda}f^{\lambda}_{\mu\nu}P_{\lambda}(x;q,t)
\label{eq:f-LR}
\end{equation}
Now let $\lambda,\mu$ be partitions and define $Q_{\lambda/\mu}\in\Lambda_{\F}$ by
\begin{equation}
Q_{\lambda/\mu}(x;q,t)
=\sum_{\nu}f^{\lambda}_{\mu\nu}Q_{\nu}(x;q,t).
\label{eq:Q-skew-def}
\end{equation}
Then $Q_{\lambda/\mu}(x;q,t)=0$ unless $\lambda\supset\mu$,
and $Q_{\lambda/\mu}$ is homogeneous of degree
$|\lambda|-|\mu|$,
which is called \defterm{Macdonald's skew $Q$-function}.
We define \defterm{Macdonald's skew $P$-function} $P_{\lambda/\mu}$ as
\begin{equation}
Q_{\lambda/\mu}(x;q,t)=\frac{b_{\lambda}(q,t)}{b_{\lambda}(q,t)}P_{\lambda/\mu}(x;q,t).
\label{eq:Q=bP-skew}
\end{equation}
holds.
Let $T$ be a tableau of shape $\lambda-\mu$ and weight $\nu$,
thought as a sequence of partitions 
$(\lambda^{(0)},\dots,\lambda^{(r)})$
such that
\[
\mu=\lambda^{(0)}\subset\lambda^{(1)}\subset\cdots\subset\lambda^{(r)}=\lambda
\]
and such that each $\lambda^{(i)}-\lambda^{(i-1)}$ is a horizontal strip.
Let
\begin{align*}
&
\phi_{T}(q,t)=\prod_{i=1}^{r}\phi_{\lambda^{(i)}/\lambda^{(i-1)}}(q,t),
\\&
\psi_{T}(q,t)=\prod_{i=1}^{r}\psi_{\lambda^{(i)}/\lambda^{(i-1)}}(q,t).
\end{align*}
Then we have \cite[VI, (7.13), (7.13')]{Mac}
\begin{align*}
&
Q_{\lambda/\mu}(x;q,t)=\sum_{T}\phi_{T}(q,t)x^{T},
\\&
P_{\lambda/\mu}(x;q,t)=\sum_{T}\psi_{T}(q,t)x^{T},
\end{align*}
summed over tableaux $T$ of shape $\lambda-\mu$,
where $x^{T}=\prod_{i=1}^{r}x_{i}^{|\lambda^{(i)}-\lambda^{(i-1)}|}$.
It also holds
\cite[VI.7, (7.9) (7.9')]{Mac}
\begin{align}
&
Q_{\lambda}(x,z;q,t)
=\sum_{\mu}Q_{\lambda/\mu}(x,z;q,t)Q_{\mu}(x,z;q,t),
\label{eq:Q-transitive}
\\&
P_{\lambda}(x,z;q,t)
=\sum_{\mu}P_{\lambda/\mu}(x,z;q,t)P_{\mu}(x,z;q,t),
\label{eq:P-transitive}
\end{align}
where the sums on the right are over partitions $\mu\subset\lambda$.
The following lemma has appeared in the proof of \cite[Proposition~2.2]{Vul2}
(also see \cite[I.5, Ex.26]{Mac} and \cite[Proposition~5.1]{Vul1}).
%
%
%
%
\begin{lemma}
\label{lem:QP}
Let $\mu$ and $\nu$ be partitions,
and $x=(x_{1},x_{2},\dots)$ and $y=(y_{1},y_{2},\dots)$
are independent indeterminates.
\begin{equation}
\sum_{\lambda}
Q_{\lambda/\mu}(x;q,t)P_{\lambda/\nu}(y;q,t)
=\Pi(x;y;q,t)
\sum_{\tau}
Q_{\nu/\tau}(x;q,t)P_{\mu/\tau}(y;q,t)
\label{eq:lem-QP}
\end{equation}
\end{lemma}
\begin{demo}{Proof}
First if we use (\ref{eq:Q-transitive}) (\ref{eq:P-transitive}) and (\ref{eq:normalisation}),
then we have
\begin{align*}
&
\sum_{\mu,\nu}\sum_{\lambda}
Q_{\lambda/\mu}(x)P_{\lambda/\nu}(y)
Q_{\mu}(z)P_{\nu}(w)
\\&=
\sum_{\mu,\nu}\sum_{\lambda}
Q_{\lambda}(x,z)P_{\lambda}(y,w)
\\&=
\Pi(x,z;y,w)
\\&=
\Pi(x;y)
\Pi(x;w)
\Pi(z;y)
\Pi(z;w)
\\&=
\Pi(x;y)
\sum_{\xi}Q_{\xi}(x)P_{\xi}(w)
\sum_{\eta}Q_{\eta}(z)P_{\eta}(y)
\sum_{\tau}Q_{\tau}(z)P_{\tau}(w)
\intertext{by (\ref{eq:f-LR}) and (\ref{eq:Q=bP})}
&=
\Pi(x;y)
\sum_{\xi,\eta,\tau}Q_{\xi}(x)P_{\eta}(y)
\sum_{\mu}\frac{b_{\eta}b_{\tau}}{b_{\mu}}f^{\mu}_{\eta\tau}Q_{\mu}(z)
\sum_{\nu}f^{\nu}_{\xi\tau}P_{\nu}(w)
\intertext{by (\ref{eq:Q-skew-def}) and (\ref{eq:Q=bP-skew})}
&=
\Pi(x;y)
\sum_{\mu,\nu,\tau}Q_{\nu/\tau}(x)P_{\mu/\tau}(y)
Q_{\mu}(z)
P_{\nu}(w).
\end{align*}
Hence, by comparing the coefficients of 
$
Q_{\mu}(z)
P_{\nu}(w)
$
in the both sides,
we obtain the desired identity.
This completes the proof.
\end{demo}
In \cite{Vul2}
Vuleti\'c has presented so-called a generalized MacMahon's formula.
The following theorem gives a generalized form of \cite[Proposition~2.2]{Vul2},
which we use in the proof of Okada's conjecture.
%
%
%
%
\begin{theorem}
\label{th:gMacMahon}
Fix a positive integer $T$ and
two partitions $\mu^{0}$ and $\mu^{T}$.
Let $x^{0},\dots,x^{T-1}$, $y^{1},\dots,y^{T}$ be sets of variables.
Then we have
\begin{align}
&\sum_{(\lambda^{1},\mu^{1},\lambda^{2},\dots,\lambda^{T})}
\prod_{i=1}^{T}Q_{\lambda^{i}/\mu^{i-1}}(x^{i-1};q,t)P_{\lambda^{i}/\mu^{i}}(y^{i};q,t)
\nonumber\\&
=\prod_{0\leq i<j\leq T}\Pi(x^{i};y^{j};q,t)
\sum_{\nu}{\textstyle Q_{\mu^{T}/\nu}(x^{0},\dots,x^{T-1};q,t)P_{\mu^{0}/\nu}(y^{1},\dots,y^{T};q,t)}
\label{eq:gMacMahon}
\end{align}
where the sum runs over $(2T-1)$-tuples 
$(\lambda^{1},\mu^{1},\lambda^{2},\dots,\mu^{T-1},\lambda^{T})$
 of partitions
satisfying
\begin{equation}
\mu^{0}\subset\lambda^{1}\supset\mu^{1}\subset\lambda^{2}\supset\mu^{2}\subset\cdots
\supset\mu^{T-1}\subset\lambda^{T}\supset\mu^{T}.
\label{eq:updown}
\end{equation}
\end{theorem}
\begin{demo}{Proof}
Our proof is similar to that of \cite[Proposition~2.2]{Vul2}.
We proceed by induction on $T$.
If $T=1$ then Lemma~\ref{lem:QP} is nothing but the desired identity.
Assume $T>1$ and (\ref{eq:gMacMahon}) holds up to $T-1$.
We need consider the sum
\begin{equation*}
S:=\sum_{(\lambda^{1},\mu^{1},\lambda^{2},\dots,\lambda^{T})}
\prod_{i=1}^{T}Q_{\lambda^{i}/\mu^{i-1}}(x^{i-1})P_{\lambda^{i}/\mu^{i}}(y^{i}),
\end{equation*}
where the sum runs over $(\lambda^{1},\mu^{1},\lambda^{2},\dots,\lambda^{T})$ satisfying
(\ref{eq:updown}).
First fix $(\mu^{1},\dots,\mu^{T-1})$ and take the sum over
$(\lambda^{1},\dots,\lambda^{T})$ using (\ref{eq:lem-QP}).
Then we obtain
\begin{equation*}
S=\prod_{i=1}^{T}\Pi(x^{i-1};y^{i})
\sum_{(\tau^{1},\mu^{1},\dots,\tau^{T})}
\prod_{i=1}^{T}Q_{\mu^{i}/\tau^{i}}(x^{i-1})P_{\mu^{i-1}/\tau^{i}}(y^{i}),
\end{equation*}
where the sum runs over $(\lambda^{1},\mu^{1},\lambda^{2},\dots,\lambda^{T})$ satisfying
\begin{equation*}
\mu^{0}\supset\tau^{1}\subset\mu^{1}\supset\tau^{2}\subset\mu^{2}\supset\cdots
\subset\mu^{T-1}\supset\tau^{T}\subset\mu^{T}.
\end{equation*}
By the induction hypothesis we can suppose
\begin{align*}
&\sum_{(\mu^{1},\tau^{2},\dots,\mu^{T-1})}
\sum\prod_{i=1}^{T}Q_{\mu^{i}/\tau^{i}}(x^{i-1})P_{\mu^{i-1}/\tau^{i}}(y^{i})
\\&
=\prod_{0\leq i<j\leq T-1}\Pi(x^{i};y^{j+1})
\sum_{\nu}{\textstyle Q_{\tau^{T}/\nu}(x^{0},\dots,x^{T-2})P_{\tau^{1}/\nu}(y^{2},\dots,y^{T})}.
\end{align*}
Hence, substituting this identity into the above $S$, we obtain
\begin{equation*}
S=\prod_{0\leq i<j\leq T}\Pi(x^{i};y^{j})
\sum_{(\tau^{1},\nu,\tau^{T})}
Q_{\mu^{T}/\tau^{T}}(x^{T-1})P_{\mu^{0}/\tau^{1}}(y^{1})
Q_{\tau^{T}/\nu}(x^{0},\dots,x^{T-2})P_{\tau^{1}/\nu}(y^{2},\dots,y^{T}),
\end{equation*}
where the sum runs over $(\tau^{1},\nu,\tau^{T})$ such that
\[
\mu^{0}\supset\tau^{1}\supset\nu\subset\tau^{T}\subset\mu^{T}.
\]
Applying (\ref{eq:Q-transitive}) and (\ref{eq:P-transitive}),
we obtain the desired identity for $T$.
This completes the theorem. 
\end{demo}
We define $P_{[\lambda,\mu]}^{\delta}(x;q,t)$ and
$Q_{[\lambda,\mu]}^{\delta}(x;q,t)$ for a pair $(\lambda,\mu)$
of partitions, a set $x=(x_1,x_2,\dots)$
of independent variables and $\delta=\pm1$ by
\[
P_{[\lambda,\mu]}^{\delta}(x;q,t)=\begin{cases}
P_{\lambda/\mu}(x;q,t)&\text{ if $\delta=+1$,}\\
Q_{\mu/\lambda}(x;q,t)&\text{ if $\delta=-1$,}
\end{cases}
\quad
Q_{[\lambda,\mu]}^{\delta}(q,t)=\begin{cases}
Q_{\lambda/\mu}(x;q,t)&\text{ if $\delta=+1$,}\\
P_{\mu/\lambda}(x;q,t)&\text{ if $\delta=-1$.}
\end{cases}
\]
Here we assume $\lambda\supset\mu$ if $\delta=+1$,
and  $\lambda\subset\mu$ if $\delta=-1$.
%
%
%
%
\begin{corollary}
\label{cor:partition-sum}
Let $n$ be a positive integer,
and $\epsilon=(\epsilon_{1},\dots,\epsilon_{n})$ a sequence of $\pm1$,
Fix a positive integer $T$ and
two partitions $\lambda^{0}$ and $\lambda^{n}$.
Let $x^{1},\dots,x^{n}$ be sets of variables.
Then we have
\begin{align}
&\sum_{(\lambda^{1},\lambda^{2},\dots,\lambda^{n-1})}
\prod_{i=1}^{n}P^{\epsilon_{i}}_{[\lambda^{i-1},\lambda^{i}]}(x^{i};q,t)
\nonumber\\&
=
\prod_{{i<j}\atop{(\epsilon_{i},\epsilon_{j})=(-1,+1)}}\Pi(x^{i};x^{j};q,t)
\sum_{\nu}
Q_{\lambda^{n}/\nu}(\{x^{i}\}_{\epsilon_{i}=-1};q,t)P_{\lambda^{0}/\nu}(\{x^{i}\}_{\epsilon_{i}=+1};q,t),
\label{eq:partition-sum1}
\\&
\sum_{(\lambda^{1},\lambda^{2},\dots,\lambda^{n-1})}
\prod_{i=1}^{n}Q^{\epsilon_{i}}_{[\lambda^{i-1},\lambda^{i}]}(x^{i};q,t)
\nonumber\\&
=
\prod_{{i<j}\atop{(\epsilon_{i},\epsilon_{j})=(-1,+1)}}\Pi(x^{i};x^{j};q,t)
\sum_{\nu}
P_{\lambda^{n}/\nu}(\{x^{i}\}_{\epsilon_{i}=-1};q,t)Q_{\lambda^{0}/\nu}(\{x^{i}\}_{\epsilon_{i}=+1};q,t),
\label{eq:partition-sum2}
\end{align}
where the sum runs over $(n-1)$-tuples 
$(\lambda^{1},\lambda^{2},\dots,\lambda^{n-1})$
 of partitions
satisfying
\begin{equation}
\begin{cases}
\lambda^{i-1}\supset\lambda^{i}&\text{ if $\epsilon_{i}=+1$,}\\
\lambda^{i-1}\subset\lambda^{i}&\text{ if $\epsilon_{i}=-1$.}
\end{cases}
\label{eq:condition-sum}
\end{equation}
\end{corollary}
%
%
\begin{demo}{Proof}
Taking $T=n$ in Theorem~\ref{th:gMacMahon},
we have
\begin{align*}
&\sum_{(\Lambda^{1},\lambda^{1},\Lambda^{2},\dots,\Lambda^{n})}
\prod_{i=1}^{n}Q_{\Lambda^{i}/\lambda^{i-1}}(X^{i-1};q,t)P_{\Lambda^{i}/\lambda^{i}}(Y^{i};q,t)
\\&
=\prod_{0\leq i<j\leq n}\Pi(X^{i};Y^{j};q,t)
\sum_{\nu}
Q_{\lambda^{n}/\nu}(X^{0},\dots,X^{n-1};q,t)P_{\lambda^{0}/\nu}(Y^{1},\dots,Y^{n};q,t)
\end{align*}
where the sum runs over 
\begin{equation*}
\lambda^{0}\subset\Lambda^{1}\supset\lambda^{1}\subset\Lambda^{2}\supset\lambda^{2}\subset\cdots
\supset\lambda^{n-1}\subset\Lambda^{n}\supset\lambda^{n}.
\end{equation*}
Put $X^{i-1}=0$ and $Y^{i}=x^{i}$ if $\epsilon_{i}=+1$,
and $X^{i-1}=x^{i}$ and $Y^{i}=0$ if $\epsilon_{i}=-1$.
Since 
$
P_{\lambda/\mu}(0;q,t)=Q_{\lambda/\mu}(0;q,t)=\begin{cases}
1&\text{ if $\lambda=\mu$,}\\
0&\text{ otherwise,}
\end{cases}
$
we obtain 
$\Lambda^{i}=\begin{cases}
\lambda^{i-1}&\text{ if $\epsilon_{i}=+1$,}\\
\lambda^{i}&\text{ if $\epsilon_{i}=-1$.}
\end{cases}
$
Hence we have the condition (\ref{eq:condition-sum}) for the sum.
Since
\[
Q_{\Lambda^{i}/\Lambda^{i-1}}(X^{i-1};q,t)P_{\Lambda^{i}/\Lambda^{i}}(Y^{i};q,t)
=\begin{cases}
P_{\lambda^{i-1}/\lambda^{i}}(x^{i};q,t)&\text{ if $\epsilon_{i}=+1$,}\\
Q_{\lambda^{i}/\lambda^{i-1}}(x^{i};q,t)&\text{ if $\epsilon_{i}=-1$.}
\end{cases}
\]
the left-hand side equals
\begin{align*}
\sum_{(\lambda^{1},\lambda^{2},\dots,\lambda^{n-1})}
\prod_{i=1}^{n}P^{\epsilon_{i}}_{[\lambda^{i-1},\lambda^{i}]}(x^{i};q,t).
\end{align*}
Meanwhile the right-hand side becomes
\[
\prod_{{i<j}\atop{(\epsilon_{i},\epsilon_{j})=(-1,+1)}}\Pi(x^{i};x^{j};q,t)
\sum_{\nu}
Q_{\lambda^{n}/\nu}(\{x^{i}\}_{\epsilon_{i}=-1};q,t)P_{\lambda^{0}/\nu}(\{x^{i}\}_{\epsilon_{i}=+1};q,t).
\]
This proves (\ref{eq:partition-sum1}).
The other identity can be proven similarly.
\end{demo}
%
%
%
%
\begin{theorem}{(Warnaar \cite[Proposition~1.3, (1.17)]{War})}
\begin{align}
\sum_{\lambda}w^{r(\lambda)}b_{\lambda}^\text{oa}(q,t)P_{\lambda}(x;q,t)
=\prod_{i\geq1}\frac{(1+w x_{i})(qtx_{i}^2;q^2)_{\infty}}{(x_{i}^2;q^2)_{\infty}}
\prod_{i<j}\frac{(tx_{i}x_{j};q)_{\infty}}{(x_{i}x_{j};q)_{\infty}},
\label{eq:Warnaar}
\end{align}
where $r(\lambda)$ is the number of rows of odd length.
\end{theorem}
Applying $w_{q,t}$ \cite[VI.2, (2.14)]{Mac} to the both sides of (\ref{eq:Warnaar}),
we obtain
%
%
\begin{corollary}
\label{cor:Cor-Warnaar}
\begin{align}
\sum_{\lambda}w^{r(\lambda')}b_{\lambda}^\text{el}(q,t)P_{\lambda}(x;q,t)
=\prod_{i\geq1}\frac{(twx_{i};q)_{\infty}}{(wx_{i};q)_{\infty}}
\prod_{i<j}\frac{(tx_{i}x_{j};q)_{\infty}}{(x_{i}x_{j};q)_{\infty}}.
\label{eq:Cor-Warnaar}
\end{align}
\end{corollary}
%
\begin{demo}{Proof}
First, if we take logarithm of the right-hand side of (\ref{eq:Warnaar}),
then we have
\begin{align*}
\log\prod_{i<j}
\frac{(tx_{i}x_{j};q)_{\infty}}{(x_{i}x_{j};q)_{\infty}}
&=\sum_{i<j}\sum_{r\geq0}
\left\{
\log\left(1-q^{r}tx_{i}x_{j}\right)
-\log\left(1-q^{r}x_{i}x_{j}\right)
\right\}
\\&=
\sum_{n\geq1}\frac1{n}\cdot\frac{1-t^{n}}{1-q^{n}}\sum_{i<j}x_{i}^{n}x_{j}^{n}
\\&=
\frac12\sum_{n\geq1}\frac1{n}\cdot\frac{1-t^{n}}{1-q^{n}}\left\{p_{n}(x)^{2}-p_{2n}(x)\right\}.
\end{align*}
Applying the $\F$-algebra homomorphism $w_{q,t}$ to this formula,
and using 
$
w_{q,t}p_{r}(x)=(-1)^{r-1}\frac{1-q^{r}}{1-t^{r}}p_{r}(x)
$
(see \cite[VI.2, (2.14)]{Mac}),
we obtain
\begin{align*}
w_{q,t}\log\prod_{i<j}
\frac{(tx_{i}x_{j};q)_{\infty}}{(x_{i}x_{j};q)_{\infty}}
&=
\sum_{n\geq1}\frac1{n}\cdot\frac{1-q^{n}t^{n}}{1-t^{2n}}\sum_{i}x_{i}^{2n}
+\sum_{n\geq1}\frac1{n}\cdot\frac{1-q^{n}}{1-t^{n}}\sum_{i<j}x_{i}^{n}x_{j}^{n}.
\end{align*}
Similarly, since
\begin{align*}
\log\prod_{i\geq1}\frac{(qtx_{i}^2;q^2)_{\infty}}{(x_{i}^2;q^2)_{\infty}}
=
\sum_{n\geq1}\frac1{n}\cdot\frac{1-q^nt^n}{1-q^{2n}}p_{2n}(x),
\end{align*}
we have
\begin{align*}
w_{q,t}
\log\prod_{i\geq1}\frac{(qtx_{i}^2;q^2)_{\infty}}{(x_{i}^2;q^2)_{\infty}}
=-\sum_{n\geq1}\frac1{n}\cdot\frac{1-q^{n}t^{n}}{1-t^{2n}}\sum_{i}x_{i}^{2n}.
\end{align*}
Finally, from
\[
\log\prod_{i}\left(1+wx_{i}\right)
=\sum_{n\geq1}\frac{(-1)^{n-1}w^{n}}{n}p_{n}(x)
\]
we obtain
\begin{align*}
w_{q,t}
\log\prod_{i}\left(1+wx_{i}\right)
&=\sum_{i}\sum_{n\geq1}\frac1{n}\cdot\frac{1-q^{n}}{1-t^{n}}w^{n}x_{i}^{n}.
\end{align*}
Hence we obtain
\[
w_{q,t}
\prod_{i\geq1}
\frac{\left(1+wx_{i}\right)\left(qtx_{i}^2;q^2\right)_{\infty}}
{\left(x_{i}^2;q^2\right)_{\infty}}
\prod_{i<j}
\frac{\left(tx_{i}x_{j};q\right)}{\left(x_{i}x_{j};q\right)}
=\prod_{i}
\frac{\left(qwx_{i};t\right)_{\infty}}{\left(wx_{i};t\right)_{\infty}}
\prod_{i<j}
\frac{\left(qx_{i}x_{j};t\right)_{\infty}}{\left(x_{i}x_{j};t\right)_{\infty}}.
\]
Now,
applying
\begin{align*}
&w_{q,t}
P_{\lambda}\left(x;q,t\right)
=Q_{\lambda'}\left(x;t,q\right)
\\
&w_{q,t}
Q_{\lambda}\left(x;q,t\right)
=P_{\lambda'}\left(x;t,q\right)
\end{align*}
(\cite[VI.5, (5.1)]{Mac})
to the left-hand side of (\ref{eq:Warnaar}),
and swapping $q$ and $t$,
we obtain the desired formula (\ref{eq:Cor-Warnaar}).
\end{demo}
From \eqref{eq:Cor-Warnaar},
we easily obtain
\begin{align}
\sum_{\lambda}w^{\frac{|\lambda|+r(\lambda')}2}b_{\lambda}^\text{el}(q,t)P_{\lambda}(x;q,t)
=\prod_{i\geq1}\frac{(twx_{i};q)_{\infty}}{(wx_{i};q)_{\infty}}
\prod_{i<j}\frac{(twx_{i}x_{j};q)_{\infty}}{(wx_{i}x_{j};q)_{\infty}},
\label{eq:Cor-Warnaar-odd}
\end{align}
and
\begin{align}
\sum_{\lambda}w^{\frac{|\lambda|-r(\lambda')}2}b_{\lambda}^\text{el}(q,t)P_{\lambda}(x;q,t)
=\prod_{i\geq1}\frac{(tx_{i};q)_{\infty}}{(x_{i};q)_{\infty}}
\prod_{i<j}\frac{(twx_{i}x_{j};q)_{\infty}}{(wx_{i}x_{j};q)_{\infty}}.
\label{eq:Cor-Warnaar-even}
\end{align}
%
%
%
%
%
%
%
%
%
%
%
%
%
%
%
%
%
%
%
%
%
%
%
%
%
%
%
%

%% file: hook03.tex
%
%
%
\section{$(q,t)$-hook formula and Macdonald polynomials}
\label{sec:proof-theorems}
In this section we rewrite the left-hand side and the right-hand side of
Okada's conjecture using the Macdonald polynomials.
In Proposition~\ref{pr:P-partitions} we give the left-hand sides
for birds and banners,
and in Proposition~\ref{pr:hook-length}
we give the right-hand sides.
We rewrite these formula into Theorem~\ref{th:left-hand-side}
and Theorem~\ref{th:right-hand-side}.
We will see Corollary~\ref{cor:partition-sum} and Corollary~\ref{cor:Cor-Warnaar} plays a key role in the proof.
\par\smallskip
We define $\phi_{[\lambda,\mu]}^{\delta}(q,t)$ and
$\psi_{[\lambda,\mu]}^{\delta}(q,t)$ for a pair $(\lambda,\mu)$ of partitions
and $\delta=\pm1$ by
\[
\phi_{[\lambda,\mu]}^{\delta}(q,t)=\begin{cases}
\phi_{\lambda/\mu}(q,t)&\text{ if $\delta=+1$,}\\
\psi_{\mu/\lambda}(q,t)&\text{ if $\delta=-1$,}
\end{cases}
\qquad
\psi_{[\lambda,\mu]}^{\delta}(q,t)=\begin{cases}
\psi_{\lambda/\mu}(q,t)&\text{ if $\delta=+1$,}\\
\phi_{\mu/\lambda}(q,t)&\text{ if $\delta=-1$.}
\end{cases}
\]
Here we assume $\lambda\succ\mu$ if $\delta=+1$,
and  $\lambda\prec\mu$ if $\delta=-1$.
We also write
\[
|\lambda-\mu|_{\delta}=\begin{cases}
|\lambda-\mu|&\text{ if $\delta=+1$,}\\
|\mu-\lambda|&\text{ if $\delta=-1$.}
\end{cases}
\]
Let $n$ be a positive integer.
Let $\epsilon=(\epsilon_{1},\dots,\epsilon_{n})$ be a sequence of $\pm1$.
Let $(\lambda^{0},\lambda^{1},\dots,\lambda^{n})$ be an $(n+1)$-tuple of partitions
such that $\lambda^{i-1}\succ\lambda^{i}$ if $\epsilon=+1$,
and $\lambda^{i-1}\prec\lambda^{i}$ if $\epsilon=-1$.
Then we write
\[
\phi^{\epsilon}_{[\lambda^{0},\lambda^{1},\dots,\lambda^{n}]}(q,t)
=\prod_{i=1}^{n}\phi_{[\lambda^{i-1},\lambda^{i}]}^{\epsilon_{i}}(q,t),
\qquad
\psi^{\epsilon}_{[\lambda^{0},\lambda^{1},\dots,\lambda^{n}]}(q,t)
=\prod_{i=1}^{n}\psi_{[\lambda^{i-1},\lambda^{i}]}^{\epsilon_{i}}(q,t).
\]
\par\smallskip
Let $\alpha$ be a strict partition,
and let $n$ be an integer such that $n\geq\alpha_1$.
Define a sequence $\epsilon=\epsilon_{n}(\alpha)=(\epsilon_1,\dots,\epsilon_n)$ of
$\pm1$ by putting
\[
\epsilon_k(\alpha)=\begin{cases}
+1 &\text{ if $k$ is a part of $\alpha$,}\\
-1 &\text{ if $k$ is not a part of $\alpha$.}
\end{cases}
\]
%
%
%
%
For example,
if $\alpha=(8,5,2,1)$ and $n=10$,
then we have $\epsilon=(++--+--+--)$.
Let $\pi\in\A(P)$ a $P$-partition
for the the shifted shape $P=P_{2}(\alpha)$.
For each integer $k=0,\dots,n$ 
we define the $k$th trace $\pi[k]$ to be the sequence
$(\dots,\pi_{2,k+2},\pi_{1,k+1})$ obtained by reading the $k$th diagonal from SE to NW.
Here we use the convention that $\pi[k]=\emptyset$ if $k\geq\alpha_{1}$.
For example,
if $\pi$ is the $P$-partition of shifted shape $\alpha=(8,5,2,1)$
in Figure~\ref{fig:Shiftedshapes-partition},
then we have
$\pi[0]=(\pi_{44},\pi_{33},\pi_{22},\pi_{11})$,
$\pi[1]=(\pi_{34},\pi_{23},\pi_{12})$,
$\pi[2]=(\pi_{24},\pi_{13})$,
$\pi[3]=(\pi_{25},\pi_{14})$,
$\pi[4]=(\pi_{26},\pi_{15})$,
$\pi[5]=(\pi_{16})$,
$\pi[6]=(\pi_{17})$,
$\pi[7]=(\pi_{18})$,
$\pi[8]=\pi[9]=\pi[10]=\emptyset$,
and
\[
\pi[0]\succ\pi[1]\succ\pi[2]\prec\pi[3]\prec\pi[4]\succ\pi[5]\prec\pi[6]\prec\pi[7]\succ\pi[8]\prec\pi[9]\prec\pi[10].
\]
By direct computation one can easily check
\begin{align*}
W_{P}(\pi;q,t)&=b_{\pi[0]}^\text{el}(q,t)
\psi^{\epsilon(\alpha)}_{[\pi[0],\dots,\pi[10]]}(q,t)
=b_{\pi[0]}^\text{el}
\psi_{\pi[0]/\pi[1]}
\psi_{\pi[1]/\pi[2]}
\phi_{\pi[3]/\pi[2]}
\\&\times
\phi_{\pi[4]/\pi[3]}
\psi_{\pi[4]/\pi[5]}
\phi_{\pi[6]/\pi[5]}
\phi_{\pi[7]/\pi[6]}
\psi_{\pi[7]/\pi[8]}
\phi_{\pi[9]/\pi[8]}
\phi_{\pi[10]/\pi[9]}.
\end{align*}
In the following
we write
\begin{align*}
&
\hPhi_{m}^{n}(\rho,\theta;q,t)=
\frac{f(\rho_{n};0)f(\theta_{n};n+1)}{f(\rho_{m};0)(\theta_{m};m+1)}
\,
\Phi_{m}^{n}(\rho,\theta;q,t),
\\&
\tPhi_{m}^{n}(\tx;\rho,\theta;q,t)=
\hPhi_{m}^{n}(\rho,\theta;q,t)
\prod_{i=m+1}^{n}\tx_{i}^{\rho_{i}+\theta_{i}-\rho_{i-1}-\theta_{i-1}}
\end{align*}
in short,
where $\rho=(\rho_{m},\dots,\rho_{n})$ and $\theta=(\theta_{m},\dots,\theta_{n})$
satisfy (\ref{eq:cond-rho-theata}),
and $\tx=(\tx_{m},\dots,\tx_{n})$ are indeterminates.
%
%
For example,
if $\pi=(\sigma,\tau;f)$ is the $P$-partition of the bird 
$P=P_{3}(\alpha,\beta;f)$
for $\alpha=(4,3)$, $\beta=(4,2)$ and $f=2$
(see Figure~\ref{fig:Birds-Banners-partition})
and satisfies (\ref{eq:PP-condition-Birds}),
then we have
\begin{align*}
W_{P}(\pi;q,t)&=\hPhi^{2}_{0}(\rho,\theta;q,t)
\psi^{\epsilon(\alpha)}_{[\sigma[0],\dots,\sigma[4]]}(q,t)
\phi^{\epsilon(\beta)}_{[\tau[0],\dots,\tau[4]]}(q,t).
\end{align*}
Further,
if $\pi=(\sigma;f)$ is the $P$-partition of the banner
$P=P_{5}(\alpha;f)$
for $\alpha=(10,6,3,2)$ and $f=2$
(see Figure~\ref{fig:Birds-Banners-partition})
and satisfies (\ref{eq:PP-condition-Banners}),
then we have
\begin{align*}
W_{P}(\pi;q,t)&=\hPhi^{2}_{1}(\rho,\theta;q,t)
b_{\sigma[0]}^\text{el}(q,t)
\psi^{\epsilon(\alpha)}_{[\sigma[0],\dots,\sigma[10]]}(q,t).
\end{align*}
%
%
%
%
\begin{prop}
\label{prop:weights}
\begin{enumerate}
\item[(1)]
Let $P=P_{2}(\alpha)$ be the shifted shape associated with
a strict partition $\alpha$ such that $\ell(\alpha)=r$,
and let $n$ be an integer such that $n\geq\alpha_1$.
If $\pi\in\A(P)$ is a $P$-partition 
 satisfying the condition (\ref{eq:Shifted-PP-condition}),
then we have
\begin{equation}
W_{P}(\pi;q,t)
=b^\text{el}_{\pi[0]}(q,t)\psi^{\epsilon(\alpha)}_{[\pi[0],\dots,\pi[n]]}(q,t)
=\frac{b^\text{el}_{\pi[0]}(q,t)}{b_{\pi[0]}(q,t)}
\phi^{\epsilon(\alpha)}_{[\pi[0],\dots,\pi[n]]}(q,t)
\label{eq:weights-ShiftedShapes}
\end{equation}
and
\begin{equation}
z^{\pi}=w^{\frac{|\pi[0]|-r(\pi[0]')}2}\prod_{i=1}^{n}\tz_{i}^{\,\,\epsilon_{i}(\alpha)\,|\pi[i-1]-\pi[i]|_{\epsilon_{i}(\alpha)}},
\label{eq:z-ShiftedShapes}
\end{equation}
where 
$w$ and $\tz_{i}$ ($1\leq i\leq n$)
are as in Proposition~\ref{pr:hook-length}~(1).
\item[(2)]
Let $\alpha=(\alpha_1,\alpha_2)$ 
and $\beta=(\beta_1,\beta_2)$
be strict partitions such that $\ell(\alpha)=\ell(\beta)=2$.
Let $f>0$ be a positive integer,
and set $P=P_{3}(\alpha,\beta;f)$ the bird associated with
$\alpha$, $\beta$ and $f$.
Let $m$ (resp. $n$) be a positive integer such that $m\geq\alpha_{1}$
(resp. $n\geq\beta_{1}$).
If $\pi=(\sigma,\tau;\rho,\theta)$ is a $P$-partition 
 satisfying the condition (\ref{eq:PP-condition-Birds}),
then we have
\begin{align}
W_{P}(\pi;q,t)
&=\hPhi_{0}^{f}(\rho,\theta;q,t)
\,
\psi^{\epsilon(\alpha)}_{[\sigma[0],\dots,\sigma[m]]}(q,t)
\,
\phi^{\epsilon(\beta)}_{[\tau[0],\dots,\tau[n]]}(q,t)
\label{eq:weights-Birds}
\end{align}
and
\begin{equation}
z^{\pi}=
\tx_{0}^{\rho_{0}+\theta_{0}}
\prod_{i=1}^{m}\tz_{i}^{\,\,\epsilon_{i}(\alpha)\,|\sigma[i-1]-\sigma[i]|_{\epsilon_{i}(\alpha)}}
\prod_{i=1}^{n}\ty_{i}^{\,\,\epsilon_{i}(\beta)\,|\tau[i-1]-\tau[i]|_{\epsilon_{i}(\beta)}}
\prod_{i=1}^{f}\tx_{i}^{\rho_{i}+\theta_{i}-\rho_{i-1}-\theta_{i-1}},
\label{eq:z-Birds}
\end{equation}
where 
$\tx_{i}$ ($0\leq i\leq f$), 
$\ty_{i}$ ($1\leq i\leq n$)
and 
$\tz_{i}$ ($1\leq i\leq m$)
are as in Proposition~\ref{pr:hook-length}~(2).
\item[(3)]
Let $\alpha=(\alpha_1,\alpha_2,\alpha_3,\alpha_4)$ 
be a strict partition such that $\ell(\alpha)=4$.
Let $P=P_{6}(\alpha;f)$ the banner associated with
$\alpha$ and $f$.
If $\pi=(\sigma;\rho,\theta)$ is a $P$-partition 
satisfying the condition (\ref{eq:PP-condition-Banners}),
then we have
\begin{align}
W_{P}(\pi;q,t)
&=
\hPhi_{1}^{f}(\rho,\theta;q,t)
\,
b^\text{el}_{\sigma[0]}(q,t)
\,
\psi^{\epsilon(\alpha)}_{[\sigma[0],\dots,\sigma[n]]}(q,t)
\label{eq:weights-Banners}
\end{align}
and
\begin{equation}
z^{\pi}=(\tx_{2}w)^{\sigma_{11}+\sigma_{33}}
\prod_{i=2}^{f}\tx_{i}^{\rho_{i}+\theta_{i}-\rho_{i-1}-\theta_{i-1}}
\prod_{i=1}^{n}\tz_{i}^{\,\,\epsilon_{i}(\alpha)\,|\sigma[i-1]-\sigma[i]|_{\epsilon_{i}(\alpha)}},
\label{eq:z-Banners}
\end{equation}
where $w$, $\tx_{i}$ ($1\leq i\leq f$),
 and $\tz_{i}$ ($1\leq i\leq n$)
are as in Proposition~\ref{pr:hook-length}~(3).
\end{enumerate}
\end{prop}
%
%
%
%
\begin{demo}{Proof}
\begin{enumerate}
\item[(1)]
From (\ref{eq:f-ND}) and (\ref{eq:Pierri-psi})
 we have
\begin{align*}
&
f_{\alpha}^\text{ND}(\pi;q,t)
=\begin{cases}
\prod_{1\leq i\leq j}
\frac
{f(\pi[1]_{i}-\pi[0]_{j+1};j-i)}
{f(\pi[1]_{i}-\pi[1]_{j};j-i)}
\,\prod_{i=2}^{n}\psi_{[\pi[i-1],\pi[i]]}^{\epsilon_{i}(\alpha)}(q,t)
&\text{ if $\epsilon_{1}(\alpha)=+$,}\\
\prod_{1\leq i\leq j}
\frac
{f(\pi[1]_{i}-\pi[0]_{j};j-i)}
{f(\pi[1]_{i}-\pi[1]_{j};j-i)}
\,\prod_{i=2}^{n}\psi_{[\pi[i-1],\pi[i]]}^{\epsilon_{i}(\alpha)}(q,t)
&\text{ if $\epsilon_{1}(\alpha)=-$.}
\end{cases}
\end{align*}
Similarly, from (\ref{eq:f-D}) and (\ref{eq:b-el})
 we have
\begin{align*}
&
f_{\alpha}^\text{D}(\pi;q,t)
=\begin{cases}
\prod_{1\leq i\leq j}
\frac
{f(\pi[0]_{i}-\pi[1]_{j};j-i)}
{f(\pi[0]_{i}-\pi[0]_{j+1};j-i)}
\,b^\text{el}_{\pi[0]}(q,t)
&\text{ if $\epsilon_{1}(\alpha)=+$,}\\
\prod_{1\leq i\leq j}
\frac
{f(\pi[0]_{i}-\pi[1]_{j+1};j-i)}
{f(\pi[0]_{i}-\pi[0]_{j+1};j-i)}
\,b^\text{el}_{\pi[0]}(q,t)
&\text{ if $\epsilon_{1}(\alpha)=-$.}
\end{cases}
\end{align*}
Hence we obtain (\ref{eq:weights-ShiftedShapes}) from
(\ref{eq:weight-SShapes}) since
\[
\psi_{[\pi[0],\pi[1]]}^{\epsilon_{1}(\alpha)}(q,t)
=\begin{cases}
\prod_{1\leq i\leq j}
\frac
{f(\pi[0]_{i}-\pi[1]_{j};j-i)f(\pi[1]_{i}-\pi[0]_{j+1};j-i)}
{f(\pi[1]_{i}-\pi[1]_{j};j-i)f(\pi[0]_{i}-\pi[0]_{j+1};j-i)}
&\text{ if $\epsilon_{1}(\alpha)=+$,}\\
\prod_{1\leq i\leq j}
\frac
{f(\pi[1]_{i}-\pi[0]_{j};j-i)f(\pi[0]_{i}-\pi[1]_{j+1};j-i)}
{f(\pi[1]_{i}-\pi[1]_{j};j-i)f(\pi[0]_{i}-\pi[0]_{j+1};j-i)}
&\text{ if $\epsilon_{1}(\alpha)=-$.}
\end{cases}
\]
Meanwhile, (\ref{eq:z-ShiftedShapes}) can be easily obtained from
\[
z^{\pi}=w^{\pi_{r-1,r-1}+\pi_{r-3,r-3}+\dots}\prod_{i=1}^{n}\tz_{i}^{\,\,|\pi[i-1]|-|\pi[i]|}.
\]
\item[(2)]
As in (1) we have
\begin{align*}
&
f_{\alpha}^\text{ND}(\sigma;q,t)
\\&
=\begin{cases}
f(\sigma_{12}-\sigma_{11};0)
\,\prod_{i=2}^{n}\psi_{[\sigma[i-1],\sigma[i]]}^{\epsilon_{i}(\alpha)}(q,t)
&\text{if $\epsilon_{1}(\alpha)=+$,}\\
\frac
{f(\sigma_{23}-\sigma_{22};0)f(\sigma_{23}-\sigma_{11};1)f(\sigma_{12}-\sigma_{11};0)}
{f(\sigma_{23}-\sigma_{12};1)}
\,\prod_{i=2}^{n}\psi_{[\sigma[i-1],\sigma[i]]}^{\epsilon_{i}(\alpha)}(q,t)
&\text{if $\epsilon_{1}(\alpha)=-$.}
\end{cases}
\end{align*}
From (\ref{eq:f-ND}) and (\ref{eq:Pierri-phi})
 we have
\begin{align*}
&
f_{\beta}^\text{ND}(\tau;q,t)
\\&
=\begin{cases}
f(\tau_{12}-\tau_{11};0)
\,\prod_{i=2}^{n}\phi_{[\tau[i-1],\tau[i]]}^{\epsilon_{i}(\beta)}(q,t)
&\text{if $\epsilon_{1}(\beta)=+$,}\\
\frac
{f(\tau_{23}-\tau_{22};0)f(\tau_{23}-\tau_{11};1)f(\tau_{12}-\tau_{11};0)}
{f(\tau_{23}-\tau_{12};0)}
\,\prod_{i=2}^{n}\phi_{[\tau[i-1],\tau[i]]}^{\epsilon_{i}(\beta)}(q,t)
&\text{if $\epsilon_{1}(\beta)=+$.}
\end{cases}
\end{align*}
%
%
Hence, if we use (\ref{eq:Pierri-phi}) or (\ref{eq:Pierri-psi}),
then we obtain (\ref{eq:weights-Birds}) from (\ref{eq:weight-Birds}).
On the other hand,
(\ref{eq:z-Birds}) is easily obtained from
\[
z^{\pi}=
z_{0}^{\,\sigma_{11}+\sigma_{22}}
\prod_{i=1}^{f}x_{i}^{\rho_{i}+\theta_{i}}
\prod_{i=1}^{m}\tz_{i}^{\,\,\epsilon_{i}(\alpha)\,|\sigma[i-1]-\sigma[i]|_{\epsilon_{i}(\alpha)}}
\prod_{i=1}^{n}\ty_{i}^{\,\,\epsilon_{i}(\beta)\,|\tau[i-1]-\tau[i]|_{\epsilon_{i}(\beta)}}
\]
using
$
z_{0}^{\,\sigma_{11}+\sigma_{22}}
\prod_{i=0}^{f}x_{i}^{\rho_{i}+\theta_{i}}=
(z_{0}\tx_{1})^{\rho_{0}+\theta_{0}}\prod_{i=1}^{f}\tx_{i}^{\,\,\rho_{i}+\theta_{i}-\rho_{i-1}-\theta_{i-1}}
$,
where we use the convention $\sigma_{11}=\rho_{0}$ and $\sigma_{22}=\theta_{0}$.
\item[(3)]
As in (1) we have
\[
f_{\alpha}^\text{D}(\sigma;q,t)
f_{\alpha}^\text{ND}(\sigma;q,t)
=
b^\text{el}_{\sigma[0]}(q,t)\psi^{\epsilon(\alpha)}_{[\sigma[0],\dots,\sigma[n]]}(q,t).
\]
Hence we can obtain (\ref{eq:weights-Banners}).
Meanwhile, (\ref{eq:z-Banners}) can be obtained from
\[
z^{\pi}=\prod_{i=2}^{f}x_{i}^{\rho_{i}+\theta_{i}}\cdot
\left(\frac{z_{0}}{y_{0}}\right)^{\sigma_{11}+\sigma_{33}}
\prod_{i=1}^{n}\tz_{i}^{\,\,|\pi[i-1]|-|\pi[i]|}
\]
using 
$
\prod_{i=2}^{f} x_{i}^{\rho_{i}+\theta_{i}}
=\tx_{2}^{\rho_{1}+\theta_{1}} \prod_{i=2}^{f} \tx_{i}^{\rho_{i}+\theta_{i}-\rho_{i}-1-\theta_{i-1}}
$,
where we use the convention $\rho_{1}=\sigma_{11}$ and $\theta_{1}=\pi_{33}$.
\end{enumerate}
\end{demo}
%
%
%
%
\begin{theorem}
\label{th:left-hand-side}
\begin{enumerate}
\item[(1)]
Let $P=P_{2}(\alpha)$ be the shifted shape
associated with a strict partition $\alpha$ of length $r$.
Let $n$ be an integer such that $n\geq\alpha_{1}$,
and let $\alpha^{c}$ be the strict partition 
formed by the complement of $\alpha$ in $[n]$.
Then we have
\begin{align}
\sum_{\pi\in\mathcal{A}\left(P\right)}
W_{P}\left(\pi;q,t\right) z^{\pi} 
=
\prod_{\alpha^{c}_{k}<\alpha_{l}}
F\left({\tz_{\alpha^{c}_{k}}}^{-1}\tz_{\alpha_{l}}\right)
\sum_{\lambda}
w^{\frac{|\lambda|-r(\lambda')}2}b_{\lambda}^\text{el}(q,t)P_{\lambda}\left(\tz_{\alpha_{1}}\dots, \tz_{\alpha_{r}} ; q, t\right),
\label{eq:LHS-ShiftedShapes}
\end{align}
where $w$ and $\tz_{i}$ ($i=1,\dots,n$) are as in Proposition~\ref{pr:hook-length}~(1).
\item[(2)]
Let $\alpha=(\alpha_1,\alpha_2)$ 
and $\beta=(\beta_1,\beta_2)$
be strict partitions such that $\ell(\alpha)=\ell(\beta)=2$.
Let $f>0$ be a positive integer,
and set $P=P_{3}(\alpha,\beta;f)$ to be the bird associated with
$\alpha$, $\beta$ and $f$.
Let $m$ (resp. $n$) be a positive integer such that $m\geq\alpha_{1}$
(resp. $n\geq\beta_{1}$).
If $\pi=(\sigma,\tau;\rho,\theta)$ is a $P$-partition 
 satisfying the condition (\ref{eq:PP-condition-Birds}),
then we have
\begin{align}
&
\sum_{\pi\in\A(P)}W_{P}(\pi;q,t)z^{\pi}
=
\prod_{\alpha_{i}^{c}<\alpha_{j}}
F\left(\tz_{\alpha_{i}^{c}}^{-1}\tz_{\alpha_{j}}\right)
\prod_{\beta_{i}^{c}<\beta_{j}}
F\left(\ty_{\beta_{i}^{c}}^{-1}\ty_{\beta_{j}}\right)
\nonumber\\&\times
\sum_{(\rho,\theta)}
\tPhi_{0}^{f}(\tx;\rho,\theta;q,t)
P_{(\theta_{0},\rho_{0})}(\tx_{0}\tz_{\alpha_{1}},\tx_{0}\tz_{\alpha_{2}};q,t)
Q_{(\theta_{0},\rho_{0})}(\ty_{\beta_{1}},\tz_{\beta_{2}};q,t).
\label{eq:LHS-Birds}
\end{align}
where the sum on the right-hand side is taken over all pairs 
$(\rho,\theta)$
with
$\rho=(\rho_{0},\dots,\rho_{f})$ and $\theta=(\theta_{0},\dots,\theta_{f})$
satisfying
\begin{equation}
\label{eq:rho-theta}
0\leq\rho_{f}\leq\cdots\leq\rho_{0}\leq\theta_{0}\leq\cdots\leq\theta_{f}.
\end{equation}
Here $\tx_{i}$ ($0\leq i\leq f$),
$y_{i}$ ($1\leq i\leq n$) and $z_{i}$ ($1\leq i\leq m$) 
are as in Proposition~\ref{pr:hook-length}~(2).
\item[(3)]
Let $\alpha=(\alpha_1,\alpha_2,\alpha_3,\alpha_4)$ 
be a strict partition such that $\ell(\alpha)=4$.
Let $P=P_{6}(\alpha;f)$ be the banner associated with
$\alpha$ and $f$.
If $\pi=(\sigma;\rho,\theta)$ is a $P$-partition 
satisfying the condition (\ref{eq:PP-condition-Banners}),
then we have
\begin{align}
\sum_{\pi\in\A(P)}W_{P}(\pi;q,t)
&=
\prod_{\lambda^{c}_{k}<\lambda_{l}}
F\left({\tz_{\lambda^{c}_{k}}}^{-1}\tz_{\lambda_{l}}\right)
\sum_{(\lambda,\rho,\theta)}
\tPhi_{1}^{f}(\tx;\rho,\theta;q,t)
\nonumber\\&\times
(\tx_{2}w)^{\lambda_{2}+\lambda_{4}}
b^\text{el}_{\lambda}(q,t)
P_{\lambda}(\tz_{\alpha_{1}},\tz_{\alpha_{2}},\tz_{\alpha_{3}},\tz_{\alpha_{4}};q,t),
\label{eq:LHS-Banners}
\end{align}
where 
 the sum on the right-hand side is taken over all triplets 
$(\lambda,\rho,\theta)$
with
$\lambda=(\lambda_{1},\lambda_{2},\lambda_{3},\lambda_{4})$,
$\rho=(\rho_{1},\dots,\rho_{f})$
and $\theta=(\theta_{1},\dots,\theta_{f})$
satisfying
\begin{equation}
\label{eq:lambda-rho-theta}
\begin{array}{lll}
\lambda_{4}\leq\lambda_{3}\leq\lambda_{2}\leq\lambda_{1},
&
0\leq\rho_{f}\leq\cdots\leq\rho_{1}=\lambda_{4},
&
\lambda_{2}=\theta_{1}\leq\cdots\leq\theta_{f}.
\end{array}
\end{equation}
Here $w$, $\tx_{i}$ ($1\leq i\leq f$) and
$z_{i}$ ($1\leq i\leq n$) 
are as in Proposition~\ref{pr:hook-length}~(3).
\end{enumerate}
\end{theorem}
\begin{demo}{Proof}
\begin{enumerate}
\item[(1)]
Since
\begin{align*}
&
\psi_{\pi[i-1]/\pi[i]}(q,t)\tz_{i}^{\,\,|\pi[i-1]-\pi[i]|}
=P_{\pi[i-1]/\pi[i]}(\tz_{i};q,t),
\\&
\phi_{\pi[i]/\pi[i-1]}(q,t)\tz_{i}^{\,\,-|\pi[i-1]-\pi[i]|}
=Q_{\pi[i]/\pi[i-1]}(\tz_{i}^{\,-1};q,t)
\end{align*}
(see \cite[VI.7, (7.14)(7.14')]{Mac}),
we can use (\ref{eq:partition-sum1})
to take the sum of the product of (\ref{eq:weights-ShiftedShapes})
and (\ref{eq:z-ShiftedShapes}),
then we obtain
\begin{align*}
&\sum_{\pi}W_{P}(\pi;q,t)z^{\pi}
\\&
=
\prod_{\alpha^{c}_{k}<\alpha_{l}}
F\left({\tz_{\alpha^{c}_{k}}}^{-1}\tz_{\alpha_{l}}\right)
\sum_{\pi[0]}
b^\text{el}_{\pi[0]}(q,t)
w^{(|\pi[0]|-r(\pi[0]'))/2}
P_{\pi[0]}(\tz_{\alpha_{1}},\dots,\tz_{\alpha_{r}};q,t),
\end{align*}
where the sum on the right-hand side runs over all partitions $\pi[0]$.
\item[(2)]
Again, using (\ref{eq:partition-sum1})
to take the sum of the product of (\ref{eq:weights-Birds})
and (\ref{eq:z-Birds}),
we obtain
\begin{align*}
&
\sum_{\pi}W_{P}(\pi;q,t)z^{\pi}
=
\prod_{\alpha^{c}_{k}<\alpha_{l}}
F\left({\tz_{\alpha^{c}_{k}}}^{-1}\tz_{\alpha_{l}}\right)
\prod_{\beta^{c}_{k}<\beta_{l}}
F\left({\ty_{\beta^{c}_{k}}}^{-1}\ty_{\beta_{l}}\right)
\tx_{0}^{\rho_{0}+\theta_{0}}
\\&\qquad\times
\sum_{(\rho,\theta)}
\tPhi_{0}^{f}(\rho,\theta)
\tx_{0}^{\rho_{0}+\theta_{0}}
P_{\sigma[0]}(\tz_{\alpha_{1}},\tz_{\alpha_{2}};q,t)
Q_{\tau[0]}(\ty_{\beta_{1}},\ty_{\beta_{2}};q,t),
\end{align*}
where the sum on the right-hand side runs over all pairs $(\rho,\theta)$
satisfying (\ref{eq:rho-theta}) with $\sigma[0]=\tau[0]=(\theta_{0},\rho_{0})$.
Finally we use 
$
\tx_{0}^{\rho_{0}+\theta_{0}}
P_{(\theta_{0},\rho_{0})}(\tz_{\alpha_{1}},\tz_{\alpha_{2}};q,t)
=
P_{(\theta_{0},\rho_{0})}(\tx_{0}\tz_{\alpha_{1}},\tx_{0}\tz_{\alpha_{2}};q,t)
$.
\item[(3)]
Using (\ref{eq:partition-sum1})
to take the sum of the product of (\ref{eq:weights-Banners})
and (\ref{eq:z-Banners}),
we obtain
\begin{align*}
&
\sum_{\pi}W_{P}(\pi;q,t)z^{\pi}
=
\prod_{\alpha^{c}_{k}<\alpha_{l}}
F\left({\tz_{\alpha^{c}_{k}}}^{-1}\tz_{\alpha_{l}}\right)
\sum_{(\rho,\theta)}
\tPhi_{1}^{f}(\sigma[0],\rho,\theta)
\\&\qquad\times
(\tx_{2}w)^{\pi[0]_{2}+\pi[0]_{4}}
b^\text{el}_{\pi[0]}(q,t)
w^{(|\sigma[0]|-r(\sigma[0]'))/2}
P_{\sigma[0]}(\tz_{\alpha_{1}},\tz_{\alpha_{2}},\tz_{\alpha_{3}},\tz_{\alpha_{4}};q,t),
\end{align*}
where the sum on the right-hand side runs over all triplets $(\sigma[0],\rho,\theta)$
satisfying (\ref{eq:lambda-rho-theta}).
\end{enumerate}
\end{demo}
If we apply Warner's formula (\ref{eq:Cor-Warnaar-even}) to (\ref{eq:LHS-ShiftedShapes}) we can obtain the $(q,t)$-hook formula (\ref{eq:product-SShapes}) for shifted shapes.
This gives another proof of \cite[Proposition~4.5~(b)]{Oka}.
Now we look at the right-hand side of the conjectured identies
in the cases of Birds and Banners.
From Proposition~\ref{pr:hook-length} we can derive the following theorem.
%
%
%
%
%
\begin{theorem}
\label{th:right-hand-side}
\begin{enumerate}
\item[(1)]
Let $\alpha=(\alpha_1,\alpha_2)$ 
and $\beta=(\beta_1,\beta_2)$
be strict partitions of length $2$.
Let $f>0$ be a positive integer,
and set $P=P_{3}(\alpha,\beta;f)$ the bird associated with
$f$, $\alpha$ and $\beta$.
Let $m,n$ be integers such that $m\geq\ell(\alpha)$ and $n\geq\ell(\beta)$,
and let $\alpha^{c}$ (resp. $\beta^{c}$) be the strict partition formed by the complement of $\alpha$ (resp. $\beta$) in $[m]$ (resp. $[n]$).
Then we have
\begin{align}
&
F\left(z[H_{p}];q,t\right)
=\prod_{\alpha_{i}^{c}<\alpha_{j}}
F\left(\tz_{\alpha_{i}^{c}}^{-1}\tz_{\alpha_{j}}\right)
\prod_{\beta_{i}^{c}<\beta_{j}}
F\left(\ty_{\beta_{i}^{c}}^{-1}\ty_{\beta_{j}}\right)
\nonumber\\&\quad\times
\sum_{{\lambda}\atop{\ell(\lambda)\leq2}}
\sum_{l=0}^{\lambda_{2}}
\sum_{k_{1},\dots,k_{f}\geq0}
\sum_{{l_{1},\dots,l_{f}\geq0}\atop{l_{1}+\dots+l_{f}=l}}
\prod_{i=1}^{f}
f(k_{i},0)f(l_{i},0)
\tx_{i}^{k_{i}-l_{i}}
\nonumber\\&\quad\times
\frac{b_{\lambda-l\cdot1^2}(q,t)}{b_{\lambda}(q,t)}
P_{\lambda}(\tx_{1}\tz_{\alpha_{1}},\tx_{1}\tz_{\alpha_{2}};q,t)
Q_{\lambda}(\ty_{\beta_{1}},\ty_{\beta_{2}};q,t)
\label{eq:RHS-Birds}
\end{align}
where $\tx_{i}$ ($1\leq i\leq f$, 
$\ty_{i}$ ($1\leq i\leq n$)
and $\tz_{i}$ ($1\leq i\leq m$) are as in Proposition~\ref{pr:hook-length} (2).
\item[(2)]
Let $\alpha=(\alpha_1,\alpha_2,\alpha_3,\alpha_4)$ 
be a strict partition of length $4$.
Let $P=P_{6}(f;\alpha)$ the Banner associated with
$\alpha$ and $\beta$.
Let $n$ be an integer such that $n\geq 4=\ell(\alpha)$,
and let $\alpha^{c}$ be the strict partition formed by the complement of $\alpha$ in $[n]$.
We write $y_{0}=z_{0'}$ and
$x_{i}=z_{-i}$ for $i=1,\dots,f$.
Then we have
\begin{align}
&F\left(z[H_{p}];q,t\right)
=\prod_{\alpha_{i}^{c}<\alpha_{j}}
F\left(\tz_{\alpha_{i}^{c}}^{-1}\tz_{\alpha_{j}}\right)
\nonumber\\&\qquad\times
\sum_{{\lambda}\atop{\ell(\lambda)\leq4}}
\sum_{l=0}^{\lambda_{4}}
\sum_{k_{2},\dots,k_{f}\geq0}
\sum_{{l_{2},\dots,l_{f}\geq0}\atop{l_{2}+\cdots+l_{f}=l}}
\prod_{i=2}^{f}
f(k_{i};0)f(l_{i};0)
\tx_{i}^{k_{i}-l_{i}}
\nonumber\\&\qquad\times
(\tx_{2}w)^{\lambda_{2}+\lambda_{4}}b_{\lambda-l\cdot1^{4}}^\text{el}(q,t)
P_{\lambda}(\tz_{\alpha_{1}},\tz_{\alpha_{2}},\tz_{\alpha_{3}},\tz_{\alpha_{4}};q,t)
\label{eq:RHS-Banners}
\end{align}
where $w$, $\tx_{i}$ ($2\leq i\leq f$) and $\tz_{i}$ ($1\leq i\leq n$)
 are as in Proposition~\ref{pr:hook-length} (3).
\end{enumerate}
\end{theorem}
\begin{demo}{Proof}
\begin{enumerate}
\item[(1)]
From (\ref{eq:normalisation}) we have
\[
\prod_{i,j=1}^{2}F\left(\tx_{1}\ty_{\beta_{j}}\tz_{\alpha_{i}}\right)
=\sum_{\mu}
P_{\mu}(\tx_{1}\tz_{\alpha_{1}},\tx_{1}\tz_{\alpha_{2}})
Q_{\mu}(\ty_{\beta_{1}},\ty_{\beta_{2}}).
\]
By the binomial theorem we have
\begin{align*}
&
\prod_{i=1}^{f}F\left(\tx_{i}\right)
=\sum_{k_{1},\dots,k_{f}\geq0}\prod_{i=1}^{f}f(k_{i};0)\tx_{i}^{k_{i}},
\\&
\prod_{i=1}^{f}F\left(\frac{\tx_{1}^{2}}{\tx_{i}}
\prod_{k,l=1}^{2}\ty_{l}\tz_{k}\right)
=\sum_{l_{1},\dots,l_{f}\geq0}\prod_{i=1}^{f}f(l_{i};0)\tx_{i}^{-l_{i}}
\left(\tx_{1}^{2}\prod_{k,l=1}^{2}\ty_{l}\tz_{k}\right)^{l_{1}+\dots+l_{f}}.
\end{align*}
By \cite[VI.4, (4.17)]{Mac} and (\ref{eq:Q=bP}) we obtain
\begin{align*}
&
\left(\tx_{1}^{2}\tz_{1}\tz_{2}\right)^{l}
P_{\mu}(\tx_{1}\tz_{\alpha_{1}},\tx_{1}\tz_{\alpha_{2}})
=
P_{\mu+l\cdot1^{2}}(\tx_{1}\tz_{\alpha_{1}},\tx_{1}\tz_{\alpha_{2}}),
\\&
\left(\ty_{1}\ty_{2}\right)^{l}
Q_{\mu}(\ty_{\beta_{1}},\ty_{\beta_{2}})
=
\frac{b_{\mu}(q,t)}{b_{\mu+l\cdot1^2}(q,t)}\,
Q_{\mu+l\cdot1^{2}}(\ty_{\beta_{1}},\ty_{\beta_{2}}).
\end{align*}
From (\ref{eq:product-Birds}) we obtain
\begin{align*}
&
F\left(z[H_{p}];q,t\right)
=\prod_{\alpha_{i}^{c}<\alpha_{j}}
F\left(\tz_{\alpha_{i}^{c}}^{-1}\tz_{\alpha_{j}}\right)
\prod_{\beta_{i}^{c}<\beta_{j}}
F\left(\ty_{\beta_{i}^{c}}^{-1}\ty_{\beta_{j}}\right)
\\&\quad\times
\sum_{l\geq0}
\sum_{{\mu}\atop{\ell(\mu)\leq2}}
\sum_{k_{1},\dots,k_{f}\geq0}
\sum_{{l_{1},\dots,l_{f}\geq0}\atop{l_{1}+\dots+l_{f}=l}}
\prod_{i=1}^{f}
f(k_{i},0)f(l_{i},0)
\tx_{i}^{k_{i}-l_{i}}
\\&\quad\times
\frac{b_{\mu}(q,t)}{b_{\mu+l\cdot1^2}(q,t)}
P_{\mu+l\cdot1^2}(\tx_{1}\tz_{\alpha_{1}},\tx_{1}\tz_{\alpha_{2}};q,t)
Q_{\mu+l\cdot1^2}(\ty_{\beta_{1}},\ty_{\beta_{2}};q,t).
\end{align*}
This immediately implies (\ref{eq:RHS-Birds}).
\item[(2)]
From Warner's formula (\ref{eq:Cor-Warnaar-even}),
we have
\begin{align*}
&
\prod_{i=1}^{4}F\left(\tz_{\alpha_{i}};q,t\right)
\prod_{1\leq i<j\leq 4}F\left(w\,\tx_{2}\,\tz_{\alpha_{i}}\tz_{\alpha_{j}};q,t\right)
\\&\qquad=
\sum_{\mu}
w^{\mu_{2}+\mu_{4}}b_{\mu}^\text{el}(q,t)
P_{\mu}(\tz_{\alpha_{1}},\tz_{\alpha_{2}},\tz_{\alpha_{3}},\tz_{\alpha_{4}};q,t).
\end{align*}
By the binomial theorem we have
\begin{align*}
&
\prod_{i=2}^{f}F\left(\tx_{i}\right)
=\sum_{k_{2},\dots,k_{f}\geq0}\prod_{i=2}^{f}f(k_{i};0)\tx_{i}^{k_{i}},
\\&
\prod_{i=2}^{f}F\left(\frac{\tx_{2}^{\,2}}{\tx_{i}}\,w^2\prod_{i=1}^{4}\tz_{\alpha_{i}}\right)
=\sum_{l_{2},\dots,l_{f}\geq0}
\prod_{i=2}^{f}f(l_{i};0)\tx_{i}^{-l_{i}}
\left(\tx_{2}^{2}w^2\prod_{i=1}^{4}\tz_{\alpha_{i}}\right)^{l_{1}+\dots+l_{f}}.
\end{align*}
From (\ref{eq:product-Birds}) we obtain
\begin{align*}
&F\left(z[H_{p}];q,t\right)
=\prod_{\alpha_{i}^{c}<\alpha_{j}}
F\left(\tz_{\alpha_{i}^{c}}^{-1}\tz_{\alpha_{j}}\right)
\\&\qquad\times
\sum_{l\geq0}
\sum_{{\mu}\atop{\ell(\mu)\leq4}}
\sum_{k_{2},\dots,k_{f}\geq0}
\sum_{{l_{2},\dots,l_{f}\geq0}\atop{l_{2}+\cdots+l_{f}=l}}
\prod_{i=2}^{f}
f(k_{i};0)f(l_{i};0)
\tx_{i}^{k_{i}-l_{i}}
\\&\qquad\times
(\tx_{2}w)^{\mu_{2}+\mu_{4}+2l}b_{\mu}^\text{el}(q,t)
P_{\mu+l\cdot1^{4}}(\tz_{\alpha_{1}},\tz_{\alpha_{2}},\tz_{\alpha_{3}},\tz_{\alpha_{4}};q,t).
\end{align*}
This immediately implies (\ref{eq:RHS-Birds}).
\end{enumerate}
\end{demo}
%
%
%
%
%
%
%
%
%
%
%
%
%
%
%
%
%
%
%

%% file: hook04.tex
%
%
%
\section{Proof by Gasper's formula}
\label{sec:proof}
Now we are in position to prove Okada's conjecture for birds and banners,
i.e.,
Theorem~\ref{th:Okada-Birds-Banners}.
At the last step of our proof
Gasper's identity (\ref{eq:Gasper}) plays an important role.
\par\smallskip
We use the fact that Macdonald's polynomials are basis for
$\Lambda_{\F}$.
(cf. \cite{LSW}).
To prove the birds case,
we fix integers $\rho_{0}$ and $\theta_{0}$
such that $\theta_{0}\geq\rho_{0}\geq0$,
and nonnegative integers $r_{1},\dots,r_{f}$.
If we compare the coefficient of
$
\prod_{i=1}^{f}\tx_{i}^{r_{i}}\cdot
P_{\lambda}(\tx_{1}\tz_{\alpha_{1}},\tx_{1}\tz_{\alpha_{2}};q,t)
Q_{\lambda}(\ty_{\beta_{1}},\ty_{\beta_{2}};q,t)
$
in (\ref{eq:LHS-Birds}) and (\ref{eq:RHS-Birds}),
the following identity must hold:
\[
\sum_{{(\rho_{1},\dots,\rho_{f})}\atop{0\leq\rho_{f}\leq\dots\leq\rho_{1}\leq\rho_{0}}}\hPhi_{0}^{f}(\rho,\theta;q,t)
=
\sum_{l=0}^{\rho_{0}}
\sum_{{l_{1},\dots,l_{f}\geq0}\atop{l_{1}+\dots+l_{f}=l}}
\frac{b_{(\theta_{0}-l,\rho_{0}-l)}(q,t)}{b_{(\theta_{0},\rho_{0})}(q,t)}
\prod_{i=1}^{f}
f(l_{i};0)f(l_{i}+r_{i};0),
\]
where $(\theta_{1},\dots,\theta_{f})$ is determined
from $\theta_{0}$ and $(\rho_{1},\dots,\rho_{f})$ by using 
the equations $\theta_{i}=\rho_{i-1}+\theta_{i-1}+r_{i}-\rho_{i}$
for $i=1,\dots,f$.
Since (\ref{eq:b}) implies
\[
b_{(\theta_{0},\rho_{0})}
=f(\theta_{0}-\rho_{0};0)\frac{f(\theta_{0};1)}{f(\theta_{0}-\rho_{0};1)}
f(\rho_{0};0),
\]
we obtain
\[
\frac{b_{(\theta_{0}-l,\rho_{0}-l)}(q,t)}{b_{(\theta_{0},\rho_{0})}(q,t)}
=\frac{f(\rho_{0}-l;0)f(\theta_{0}-l;1)}
{f(\rho_{0};0)f(\theta_{0};1)}.
\]
Hence it is enough to prove
\begin{equation}
\sum_{{(\rho_{1},\dots,\rho_{f})}\atop{0\leq\rho_{f}\leq\dots\leq\rho_{1}\leq\rho_{0}}}\hPhi_{0}^{f}(\rho,\theta;q,t)
=
\sum_{l=0}^{\rho_{0}}
\sum_{{l_{1},\dots,l_{f}\geq0}\atop{l_{1}+\dots+l_{f}=l}}
\frac{f(\rho_{0}-l;0)f(\theta_{0}-l;1)}{f(\rho_{0};0)f(\theta_{0};1)}
\prod_{i=1}^{f}
f(l_{i};0)f(l_{i}+r_{i};0).
\label{eq:Birds-Final}
\end{equation}
In the case of banners 
we fix a partition
 $\lambda=(\lambda_{1},\lambda_{2},\lambda_{3},\lambda_{4})$
of length $4$
and nonnegative integers $r_{2},\dots,r_{f}$.
If we compare the coefficient of
$
\prod_{i=2}^{f}\tx_{i}^{r_{i}}\cdot
P_{\lambda}(\tx_{1}\tz_{\alpha_{1}},\tx_{1}\tz_{\alpha_{2}};q,t)
$
in (\ref{eq:LHS-Banners}) and (\ref{eq:RHS-Banners}),
the following identity must hold:
\[
\sum_{{(\rho_{2},\dots,\rho_{f})}\atop{0\leq\rho_{f}\leq\cdots\leq\rho_{2}\leq\rho_{1}}}
\hPhi_{1}^{f}(\rho,\theta;q,t)
=
\sum_{l=0}^{\lambda_{4}}
\sum_{{l_{2},\dots,l_{f}\geq0}\atop{l_{2}+\cdots+l_{f}=l}}
\frac{b_{\lambda-l\cdot1^{4}}^\text{el}(q,t)}{b_{\lambda}^\text{el}(q,t)}
\prod_{i=2}^{f}
f(l_{i};0)f(l_{i}+r_{i};0),
\]
where $(\theta_{2},\dots,\theta_{f})$ is determined
from $\theta_{1}$ and $(\rho_{2},\dots,\rho_{f})$ by using 
the equations $\theta_{i}=\rho_{i-1}+\theta_{i-1}+r_{i}-\rho_{i}$
for $i=2,\dots,f$.
Here we use the convention that $\rho_{1}=\lambda_{4}$ and $\theta_{1}=\lambda_{2}$.
Again, because of (\ref{eq:b-el})
we obtain
\[
\frac{b_{\lambda-l\cdot1^{4}}^\text{el}(q,t)}{b_{\lambda}^\text{el}(q,t)}
=
\frac
{f(\lambda_{4}-l;0)f(\lambda_{2}-l;2)}
{f(\lambda_{4};0)f(\lambda_{2};2)}.
\]
Hence it is enough to prove
\begin{equation}
\sum_{{(\rho_{2},\dots,\rho_{f})}\atop{0\leq\rho_{f}\leq\cdots\leq\rho_{2}\leq\rho_{1}}}
\hPhi_{1}^{f}(\rho,\theta;q,t)
=
\sum_{l=0}^{\lambda_{4}}
\sum_{{l_{2},\dots,l_{f}\geq0}\atop{l_{2}+\cdots+l_{f}=l}}
\frac
{f(\lambda_{4}-l;0)f(\lambda_{2}-l;2)}
{f(\lambda_{4};0)f(\lambda_{2};2)}
\prod_{i=2}^{f}
f(l_{i};0)f(l_{i}+r_{i};0).
\label{eq:Banners-Final}
\end{equation}
In fact a more general formula holds.
If we prove the following theorem,
then the proof of (\ref{eq:Birds-Final}) and (\ref{eq:Banners-Final})
are both done.
%
%
%
%
\begin{theorem}
\label{th:lemma}
Let $m$ and $n$ be nonnegative integers.
Let $k_{0}$, $\rho_{0}$, $\theta_{0}$ be integers such that 
$0\leq k_{0}\leq\rho_{0}\leq\theta_{0}$,
and let $\gamma_{1},\dots,\gamma_{n}$ be nonnegative integers.
Then we have
\begin{align}
&
\sum_{
{(\rho_{1},\dots,\rho_{n})}
\atop
{k_{0}\leq \rho_{n}\leq \dots\leq\rho_{1}\leq\rho_{0}}
}
f(\rho_{n}-k_{0};0)
f(\theta_{n}-k_{0};m+n)
\nonumber\\&\times
\prod_{i=1}^{n}
\frac{
f(\rho_{i-1}-\rho_{i};0)
f(\theta_{i-1}-\rho_{i};i+m-1)
f(\theta_{i}-\rho_{i-1};i+m-1)
f(\theta_{i}-\theta_{i-1};0)
}
{
f(\theta_{i}-\rho_{i};i+m-1)
f(\theta_{i}-\rho_{i};i+m)
}
\nonumber\\&=
\sum_{{k_{1},\dots,k_{n}\geq0}\atop{k_{1}+\dots+k_{n}\leq\rho_{0}-\rho_{m+1}}}
f(\rho_{0}-\sum_{i=0}^{n}k_{i};0)
f(\theta_{0}-\sum_{i=0}^{n}k_{i};m)
\prod_{i=1}^{n} f(k_{i};0)f(k_{i}+\gamma_{i};0),
\label{eq:general}
\end{align}
where the sum on the left-hand side runs over all $n$-tuples 
$(\rho_{1},\dots,\rho_{n})$ of nonnegative integers such that
$k_{0}\leq \rho_{n}\leq \dots\leq\rho_{1}\leq\rho_{0}$,
 the sum on the right-hand side runs over all $n$-tuples 
$(k_{1},\dots,k_{n})$ of nonnegative integers which satisfy
$k_{1}+\dots+k_{n}\leq\rho_{0}-\rho_{m+1}$,
and
$\theta_{i}$ is determined from $\rho_{i}$, $\rho_{ii-1}$ and $\theta_{i-1}$  by
$\theta_{i}=\gamma_{i}+\theta_{i-1}+\rho_{i-1}-\rho_{i}$ for $i=1,\dots,n$.
\end{theorem}
Before we prove this theorem,
we need the following lemma which is a special case (i.e., $n=1$) of this theorem.
%
%
%
\begin{lemma}
\label{lem:special-case}
Let $m$ be a nonnegative integer.
Let $k_{0}$, $\rho_{0}$ and $\theta_{0}$ be integers such that
$0\leq k_{0}\leq\rho_{0}\leq\theta_{0}$,
and let $\gamma$ be a nonnegative integer.
Then we have
\begin{align}
&
\sum_{\rho=k_{0}}^{\rho_{0}}
f(\rho-k_{0};0)
f(\theta-k_{0};m+1)
\frac{
f(\rho_{0}-\rho;0)
f(\theta_{0}-\rho;m)
f(\theta-\rho_{0};m)
f(\theta-\theta_{0};0)
}
{
f(\theta-\rho;m)
f(\theta-\rho;m+1)
}
\nonumber\\&=
\sum_{k=0}^{\rho_{0}-k_{0}}
f(\rho_{0}-k_{0}-k;0)
f(\theta_{0}-k_{0}-k;m)
f(k;0)
f(k+\gamma;0),
\label{eq:lemma}
\end{align}
where
$
\theta=\gamma+\rho_{0}+\theta_{0}-\rho
$.
\end{lemma}
%
%
%
\begin{demo}{Proof}
Set $S_{1}$ to be the left-hand side of \eqref{eq:lemma}.
If one puts $k=\rho_{0}-\rho$,
then $\rho=\rho_{0}-k$ and $\theta=k+\gamma+\theta_{0}$.
Hence one obtains
\begin{align*}
S_{1}=&\sum_{k=0}^{\rho_{0}-k_{0}}
f(\rho_{0}-k_{0}-k;0)
f(k+\gamma+\theta_{0}-k_{0};m+1)
\\&\times
\frac{
f(k;0)
f(k+\gamma+\theta_{0}-\rho_{0};m)
f(k+\theta_{0}-\rho_{0};m)
f(k+\gamma;0)
}
{
f(2k+\gamma+\theta_{0}-\rho_{0};m)
f(2k+\gamma+\theta_{0}-\rho_{0};m+1)
}.
\end{align*}
If we use
\[
(\alpha;q)_{2k}=(\alpha^{\frac12};q)_{k}(-\alpha^{\frac12};q)_{k}(\alpha^{\frac12}q^{\frac12};q)_{k}(-\alpha^{\frac12}q^{\frac12};q)_{k},
\]
then the factors in the denominator are written as
$
f(2k+\gamma+\theta_{0}-\rho_{0};m)
=
f(\gamma+\theta_{0}-\rho_{0};m)\times\\
\frac{
\left(t^{\frac{m+1}2}q^{\frac{\gamma+\theta_{0}-\rho_{0}}2},-t^{\frac{m+1}2}q^{\frac{\gamma+\theta_{0}-\rho_{0}}2},t^{\frac{m+1}2}q^{\frac{\gamma+\theta_{0}-\rho_{0}+1}2},-t^{\frac{m+1}2}q^{\frac{\gamma+\theta_{0}-\rho_{0}+1}2};q\right)_{k}
}
{
\left(t^{\frac{m}2}q^{\frac{\gamma+\theta_{0}-\rho_{0}+1}2},-t^{\frac{m}2}q^{\frac{\gamma+\theta_{0}-\rho_{0}+1}2},t^{\frac{m}2}q^{\frac{\gamma+\theta_{0}-\rho_{0}+2}2},-t^{\frac{m}2}q^{\frac{\gamma+\theta_{0}-\rho_{0}+2}2};q\right)_{k}
}
$
and
$
f(2k+\gamma+\theta_{0}-\rho_{0};m+1)
=
f(\gamma+\theta_{0}-\rho_{0};m+1)\times\\
\frac{
\left(t^{\frac{m+2}2}q^{\frac{\gamma+\theta_{0}-\rho_{0}}2},-t^{\frac{m+2}2}q^{\frac{\gamma+\theta_{0}-\rho_{0}}2},t^{\frac{m+2}2}q^{\frac{\gamma+\theta_{0}-\rho_{0}+1}2},-t^{\frac{m+2}2}q^{\frac{\gamma+\theta_{0}-\rho_{0}+1}2};q\right)_{k}
}
{
\left(t^{\frac{m+1}2}q^{\frac{\gamma+\theta_{0}-\rho_{0}+1}2},-t^{\frac{m+1}2}q^{\frac{\gamma+\theta_{0}-\rho_{0}+1}2},t^{\frac{m+1}2}q^{\frac{\gamma+\theta_{0}-\rho_{0}+2}2},-t^{\frac{m+1}2}q^{\frac{\gamma+\theta_{0}-\rho_{0}+2}2};q\right)_{k}
}
$.
Meanwhile, the factors in the numerator are
$
f(\rho_{0}-k_{0}-k;0)
=
f(\rho_{0}-k_{0};0)
\frac{(q^{-\rho_{0}+k_{0}};q)_{k}}{(t^{-1}q^{-\rho_{0}+k_{0}+1};q)_{k}}
\left(\frac{q}{t}\right)^{k}
$,
$
f(k+\gamma+\theta_{0}-k_{0};m+1)
=
f(\gamma+\theta_{0}-k_{0};m+1)
\frac{(t^{m+2}q^{\gamma+\theta_{0}-k_{0}};q)_{k}}
{(t^{m+1}q^{\gamma+\theta_{0}-k_{0}+1};q)_{k}}
$,
$
f(k+\gamma+\theta_{0}-\rho_{0};m)
=
f(\gamma+\theta_{0}-\rho_{0};m)
\frac{(t^{m+1}q^{\gamma+\theta_{0}-\rho_{0}};q)_{k}}
{(t^{m}q^{\gamma+\theta_{0}-\rho_{0}+1};q)_{k}}
$,
$
f(k+\theta_{0}-\rho_{0};m)
=
f(\theta_{0}-\rho_{0};m)
\frac{(t^{m+1}q^{\theta_{0}-\rho_{0}};q)_{k}}
{(t^{m}q^{\theta_{0}-\rho_{0}+1};q)_{k}}
$,
$
f(k+\gamma;0)
=
f(k+\gamma;0)
\frac{(tq^{\gamma};q)_{k}}{(q^{\gamma};q)_{k}}
$.
Hence, substituting these factors, we obtain
\begin{align*}
S_{1}=&C\cdot
{}_{12}W_{11}\Bigl(bc/d;(bcq/ad)^{\frac12},-(bcq/ad)^{\frac12},
q(bc/d)^{\frac12},-q(bc/d)^{\frac12},
\\&\qquad\qquad
 ab/d,ac/d,a,b,c;q,q/a\Bigr),
\end{align*}
where $a=t$, $b=tq^{\gamma}$, $c=q^{-\rho_{0}+k_{0}}$, 
$d=t^{-m}q^{-\theta_{0}+k_{0}}$
and
\begin{align*}
C&=
\frac{
f(\rho_{0}-k_{0};0)
f(\gamma+\theta_{0}-k_{0};m+1)
f(\theta_{0}-\rho_{0};m)
f(\gamma;0)
}
{
f(\gamma+\theta_{0}-\rho_{0},m+1)
}.
\end{align*}
On the other hand,
Set $S_{2}$ to be the right-hand side  of \eqref{eq:lemma}.
If we use
$
f(\rho_{0}-k_{0}-k;0)
=
f(\rho_{0}-k_{0};0)
\frac{(q^{-\rho_{0}+k_{0}};q)_{k}}{(t^{-1}q^{-\rho_{0}+k_{0}+1};q)_{k}}
\left(\frac{q}{t}\right)^{k}
$,
$
f(\theta_{0}-k_{0}-k;m)
=
f(\theta_{0}-k_{0};m)
\frac{(t^{-m}q^{-\theta_{0}+k_{0}};q)_{k}}{(t^{-m-1}q^{-\theta_{0}+k_{0}+1};q)_{k}}
\left(\frac{q}{t}\right)^{k}
$
and
$
f(k+\gamma;0)
=
f(k+\gamma;0)
\frac{(tq^{\gamma};q)_{k}}{(q^{\gamma};q)_{k}}
$,
then we obtain
\begin{align*}
S_{2}=f(\rho_{0}-k_{0};0)f(\theta_{0}-k_{0};m)f(\gamma;0)
{}_{4}\phi_{3}\left[{
{q^{-\rho_{0}+k_{0}},t^{-m}q^{-\theta_{0}+k_{0}},t,tq^{\gamma}}
\atop
{t^{-1}q^{-\rho_{0}+k_{0}+1},t^{-m-1}q^{-\theta_{0}+k_{0}+1},q^{\gamma+1}}
};q,\frac{q^2}{t^2}\right].
\end{align*}
Hence Gasper's formula \eqref{eq:Gasper} proves that $S_{1}=S_{2}$.
The details are left to the reader.
This completes our proof.
\end{demo}
%
%
%
\begin{demo}{Proof of Theorem~\ref{th:lemma}}
We proceed by induction on $n$.
If $n=1$, then  (\ref{eq:general}) is nothing but (\ref{eq:lemma}).
Let $n\geq2$ and assume (\ref{eq:general}) is true for $n-1$.
If we set $S$ to be the left-hand side of (\ref{eq:general}),
then we have
\begin{align*}
S=
&\sum_{\rho_{1}=k_{0}}^{\rho_{0}}
\frac{
f(\rho_{0}-\rho_{1};0)
f(\theta_{0}-\rho_{1};m)
f(\theta_{1}-\rho_{0};m)
f(\theta_{1}-\theta_{0};0)
}
{
f(\theta_{1}-\rho_{1};m)
f(\theta_{1}-\rho_{1};m+1)
}
\\
&\times\sum_{
{(\rho_{2},\dots,\rho_{n})}\atop
{k_{0}\leq\rho_{n}\leq \dots\leq\rho_{2}\leq\rho_{1}}
}
f(\rho_{n}-k_{0};0)
f(\theta_{n}-k_{0};m+n)
\nonumber\\&\times
\prod_{i=2}^{n}
\frac{
f(\rho_{i-1}-\rho_{i};0)
f(\theta_{i-1}-\rho_{i};i+m-1)
f(\theta_{i}-\rho_{i-1};i+m-1)
f(\theta_{i}-\theta_{i-1};0)
}
{
f(\theta_{i}-\rho_{i};i+m-1)
f(\theta_{i}-\rho_{i};i+m)
}.
\end{align*}
We can use our induction hypothesis to obtain
\begin{align*}
S=
&
\sum_{{k_{2},\dots,k_{n}\geq0}\atop{k_{2}+\dots+k_{n}\leq\rho_{0}-k_{0}}}
\prod_{i=2}^{n} f(k_{i};0)f(k_{i}+\gamma_{i};0)
\\&\times
\sum_{\rho_{1}=k_{0}+\sum_{i=2}^{n}k_{i}}^{\rho_{0}}
f(\rho_{1}-k_{0}-\sum_{i=2}^{n}k_{i};0)
f(\theta_{1}-k_{0}-\sum_{i=2}^{n}k_{i};m+1)
\\
&\times
\frac{
f(\rho_{0}-\rho_{1};0)
f(\theta_{0}-\rho_{1};m)
f(\theta_{1}-\rho_{0};m)
f(\theta_{1}-\theta_{0};0)
}
{
f(\theta_{1}-\rho_{1};m)
f(\theta_{1}-\rho_{1};m+1)
}.
\end{align*}
If we use (\ref{eq:lemma}) again,
then we obtain
\begin{align*}
S=
&
\sum_{{k_{2},\dots,k_{n}\geq0}\atop{k_{2}+\dots+k_{n}\leq\rho_{0}-k_{0}}}
\prod_{i=2}^{n} f(k_{i};0)f(k_{i}+\gamma_{i};0)
\\&\times
\sum_{0\leq k_{1}\leq \rho_{0}-k_{0}-\sum_{i=2}^{n}k_{i}}
f(\rho_{0}-\sum_{i=0}^{n}k_{i},0)
f(\theta_{0}-\sum_{i=0}^{n}k_{i},m)
f(k_{1},0)f(k_{1}+\gamma_{1},0),
\end{align*}
which equals the right-hand side of (\ref{eq:general}).
This completes our proof.
\end{demo}
\noindent
{\bf Concluding Remarks}
In the proof of the $(q,t)$ hook formula for birds and banners,
Gasper's identity (\ref{eq:Gasper}) for ${}_{12}W_{11}$
plays an important role.
The author tried the other classes of irreducible 
$d$-complete posets, but it seems that another identity will be needed
for the rest.
\par
%
%
%
%
%
%
%
%
%
%
%